\RequirePackage{fix-cm}
\documentclass[smallcondensed]{svjour3}     \smartqed  \usepackage{graphicx}
\usepackage{soul}
\usepackage{mathptmx}      \usepackage[utf8]{inputenc}
\usepackage{amsmath}
\usepackage{bbm} \usepackage{amssymb}
\usepackage{amsfonts}
\usepackage{mathrsfs}
\usepackage{mathtools}
\usepackage{hyperref} \usepackage{ifthen} \usepackage[dvipsnames]{xcolor} \usepackage{float} \usepackage{algorithm} \usepackage{algorithmic} \usepackage{diagbox} \usepackage{pdfpages} \usepackage{stmaryrd} \usepackage{bm} \usepackage[normalem]{ulem}
\newcommand{\Rn}[1]
	{\ifthenelse{\equal{#1}{}}
		{\ensuremath{\mathbb{R}}}
		{\ensuremath{\mathbb{R}^{#1}}}}

\newcommand{\norm}[2]
	{\ifthenelse{\equal{#2}{}}
		{\ensuremath{\left\| #1 \right\|}}
		{\ensuremath{\left\| #1 \right\|_{#2}}}}

\newcommand{\argmax}{\ensuremath{\operatornamewithlimits{argmax}}}
		
\newcommand{\argmin}{\operatornamewithlimits{argmin}}

\newcommand{\rf}[1]{{\color{black} #1}}
\newcommand{\ld}[1]{{\color{black} #1}}
\newcommand{\pat}[1]{{\color{black} #1}}
\newcommand{\kl}[1]{{\color{black} #1}}

\renewcommand{\vec}[1]{{\mathbf{#1}}}
\renewcommand{\tens}[1]{{\mathbf{#1}}}

    \newcommand{\R}{\ensuremath{\mathbb{R}}}

\newcommand{\T}{\ensuremath{T}}

\newcommand{\noise}{\xi}
\newcommand{\level}{z}
\newcommand{\co}{\overline{c}}
\newcommand{\cd}{\underline{c}}

\renewcommand{\H}{\ensuremath{{\bf H}}}

\newcommand{\astar}{\ensuremath{a^*}}
\newcommand{\Sstar}{\ensuremath{S^*}}

\newcommand{\ahat}{\ensuremath{\hat{a}}}

\renewcommand{\l}{\ensuremath{\ell}}
\newcommand{\pusht}[2]{\ensuremath{#1^{\to #2}}}

\spnewtheorem{assumption}{Assumption}{\bf}{\it}

\newtheorem{cor}{Corollary}

\journalname{}
\begin{document}  

\title{Sliding window strategy for convolutional spike sorting  with Lasso}\subtitle{Algorithm, theoretical guarantees and complexity}

\titlerunning{Working set strategy for convolutional spike sorting with Lasso}        
\author{ Laurent Dragoni \and Rémi Flamary \and Karim Lounici \and Patricia Reynaud-Bouret}

\institute{Laurent Dragoni \at
            Université Côte d'Azur, CNRS, Laboratoire J.A. Dieudonné, 06108 Nice, France\\
            Tel.: +33 (0)4 89 15 62 78\\
                        \email{laurent.dragoni@univ-cotedazur.fr}                       \and
            Rémi Flamary \at
            École Polytechnique, Centre de Mathématiques Appliquées, 91128 Palaiseau, France \\
            Tel.: +33 (0)1 69 33 46 25\\
            \email{remi.flamary@polytechnique.edu}
            \and
            Karim Lounici \at
            École Polytechnique, Centre de Mathématiques Appliquées, 91128 Palaiseau, France \\
            Tel.: +33 (0)1 69 33 46 48\\
            \email{karim.lounici@polytechnique.edu}
          \and
          Patricia Reynaud-Bouret \at
            Université Côte d'Azur, CNRS, Laboratoire J.A. Dieudonné, 06108 Nice, France\\
            Tel.: +33 (0)4 89 15 04 96\\
            \email{patricia.reynaud-bouret@univ-cotedazur.fr}
          }

\date{Received: date / Accepted: date}

\maketitle

\begin{abstract}
\pat{ Spike sorting is a class of algorithms used in neuroscience to attribute the time occurences of particular electric signals, called action potential or spike, to neurons. We rephrase this problem as a particular optimization problem : Lasso for convolutional
  models in high dimension. Lasso (i.e. least absolute shrinkage and selection operator) is a very generic tool in machine learning that help us to look for sparse solutions (here the time occurrences). However, for the size of the problem at hand in this neuroscience context, the classical Lasso solvers are failing.
  We present here a new and much faster algorithm.} 
    Making use of biological properties related to
  neurons, we explain how the particular structure of the problem allows several
  optimizations, leading to an algorithm with a temporal complexity which
  grows linearly with respect to the size of the recorded signal and can
  be performed online. Moreover
  the spatial separability of the initial problem allows to break it into
  subproblems, further reducing the complexity and making possible its
  application on the latest recording devices which comprise a large number of
  sensors. We provide several mathematical results: the size and numerical
  complexity of the subproblems
  can be estimated mathematically by using percolation theory. We also show
  under reasonable assumptions that the Lasso estimator   retrieves the
  true time occurrences of the spikes {with large probability}. Finally the theoretical time complexity of the algorithm is
  given. Numerical simulations are also provided in order to illustrate the
  efficiency of our approach.


\end{abstract}

\section{Introduction}
\label{sec:introduction}

The main focus of the field of neuroscience consists in better understanding how
the brain works. From brain activity recordings, obtained for instance by
extracellular electrodes, the goal of the experimenter is to analyze the
recorded signals in order to gain insights about important aspects of the neural
code, for instance the firing rate coding or the synchronization between neurons
\cite{albert2016,lambert2018,eytan2006synchro}.
It is well established that neurons communicate by emitting shorts bursts of
current, called action potentials \pat{or spikes}. Unfortunately, in practice the recorded
signals are often a mixture of these action potentials from the neuronal
population. A fundamental pre-processing step allowing further analyses aims at
extracting from the recorded signals the mixture and spiking
activity of each neurons. This
step is called \emph{spike sorting}.
The main idea relies on a important property: since the shapes of the action
potentials emitted by the individual neurons do not vary too much along time, these
shapes can be viewed
as the signature of a particular neuron. Therefore the goal of spike sorting is
to gather the shapes of the action potentials, to associate them to their
respective neuron and to  detect at which times each particular neuron has emitted
an action potential (\pat{the time occurrences of the action potentials of a given neuron is called a spike train}).
Numerous works have focused on designing efficient spike sorting procedures lately \cite{lewicki1998review,ekanadham2011recovery,ekanadham2014unified}. The estimation of the shapes of the action potentials usually involves dimensionality reduction techniques and clustering. \pat{In simple terms, this means that the raw signal as recorded by the electrodes is cut into pieces that look like spikes. These pieces are then represented in a low dimensional space (typically via Principal Component Analysis) and grouped into homogeneous clusters, that could be understood as the spike train of a given neuron up to the noise \cite{pouzat2014}.} Despite their popularity in the neuroscience community, these methods have some severe limitations. A large number of manual calibrations from the experimenter is often required in order to recover both the shapes and the spike trains. These tasks become even harder when \pat{synchronization, that is the fact that two or more neurons emit spikes at almost the same time, takes place. But on the other hand, synchronization (see \cite{albert2016} and the references therein) is known to be fundamental for the way neurons in our brain encode information and this phenomenon cannot be discarded in the preprocessing of the data. However, in practice, synchronization in spike sorting is often solved manually later on \cite{pouzat2014}.}
Therefore, the outcomes of these methods heavily depend on the person using them \cite{wood2004variability,harris2000accuracy}. These shortcomings may also limit the usability of the most recent acquisition devices called {\it MultiElectrode Arrays} (MEA) that can contain up to several thousand electrodes, whereas the classical recordings use only 4 electrodes, with devices called {\it tetrodes}. These MEA devices tend to produce large datasets in which the spike synchronization phenomenon worsens \cite{einevoll2012mea}, making the manual sorting mentioned above infeasible.

In this paper, we propose a spike sorting procedure that is able to analyze large datasets without the need for a large number of manual calibrations. More specifically, our procedure relies on a convolutional model in order to link the activity of the neurons and the recorded signals. This model aims at reconstructing the recorded signal as a temporal convolution between the shapes of the action potentials and the sparse neurons activations.
This type of model, originally proposed in
\cite{taylor1979deconvolution} and applied for instance to speech signals
\cite{smaragdis2007}, has also been used with success  for spike sorting
\cite{ekanadham2011recovery,ekanadham2014unified}. \pat{A precursor of this model in neuroscience is due to Roberts \cite{roberts} in 1979.} A nice characteristic of this
model is based on the linearity of the convolution operator. This permits to
treat synchronizations as additive superpositions of the shapes, which is not
possible with traditional spike sorting relying on clustering
\cite{lewicki1998review}. We focus here on the estimation of the activations of
the neurons by assuming that the shapes of action potentials are known, a problem commonly refered to convolutional sparse coding.
We
explain in the next section why the hard part is indeed the estimation of the
activations. We also see that this estimation leads to a large scale Lasso convex optimization
problem, for which well-established strategies exist. \pat{Lasso stands for {\it Least Absolute Shrinkage and Selection Operator}. This learning method has been introduced in 1996 by Tibshirani \cite{tibshirani1996regression}. This method, its variants and the numerous affiliated algorithms (basis pursuit, LARS, working set, FSS, etc) have had a lot of interest in the statistical and machine learning community because of its computability and  the sparseness of its solution. In the present set-up,}
the basis pursuit estimator proposed in \cite{ekanadham2011recovery}, which is
equivalent to the Lasso estimator, verifies important statistical properties.
For example, under specific conditions \cite{BRT09,bunea2008,lounici2008}, it
can be shown that it is able to recover the support of the activations (i.e. the
activation times and the active neurons). Unfortunately computing this estimator
is hard in practice since its computational
complexity increases cubically with the number of variables of the
optimization problem, which is here proportional to the length of the signal. 
In the more general context of sparse coding problems, \cite{lee2007efficient} proposed the Feature Sign Search (FSS) algorithm, which essentially consists in a working set strategy paired with the resolution of a \pat{Quadratic Programming} problem. They notably showed that this strategy grants better results than the LARS on practical applications.
\cite{grosse2012shift} proposed an extension of FSS to convolutional sparse coding applied to audio classification, by splitting the signal into smaller temporal windows. Although they managed to improve the performances of the original FSS strategy on large signals, the fact that these windows were fixed a priori required to consider multiple passes on the whole signal.
Another similar approach, using an a priori partition of the signal proposed by \cite{moreau2018dicod}, aimed at solving smaller supbroblems with a local coordinate descent algorithm. But since the structure of this strategy imposes to wait for the convergence on every windows to reach global convergence, its domains of applications remain limited to offline analyses.

After introducing the convolutional model, we present an efficient algorithm for the computation of the Lasso. This algorithm refines the working set strategy by using temporal sliding windows, allowing it to scale in high dimension. 
In contrast with \cite{grosse2012shift,moreau2018dicod}, our method requires a single pass on the whole signal and is fast enough to allow an online analysis.

Then we explain how we can take advantage of some biological facts in
order to prove that the Lasso enjoys nice statistical properties. Furthermore we derive theoretical time complexity of the algorithm. In particular,
we note that it grows linearly with respect to the size of the recorded
signal. Other works have tackled  similar optimization problems
\cite{jas2017,la2018multivariate}, but to the best of our knowledge the approach
that we present is the first to fully take advantage of the structure of the
problem and to attain linear complexity.
Interestingly since we only perform linear operations such as convolutions, this
algorithm can be adapted to GPU architectures, allowing a very efficient scaling. Finally we
present some numerical results illustrating both the theoretical results and the
performance of our approach. Note that this algorithm is designed to solve an
estimation problem associated to a convolutional model, and it can potentially be
used to other domains, as long as the quantities of interest verifies similar
sparsity properties as in the spike sorting problem (for instance the
recognition of musical notes).

\section{Physical model \& Optimization}
\label{sec:model_optim}

\subsection{Convolutional model}
\label{subsec:convolutional_model}

\begin{figure}    \begin{center}
        \includegraphics[width=340pt]{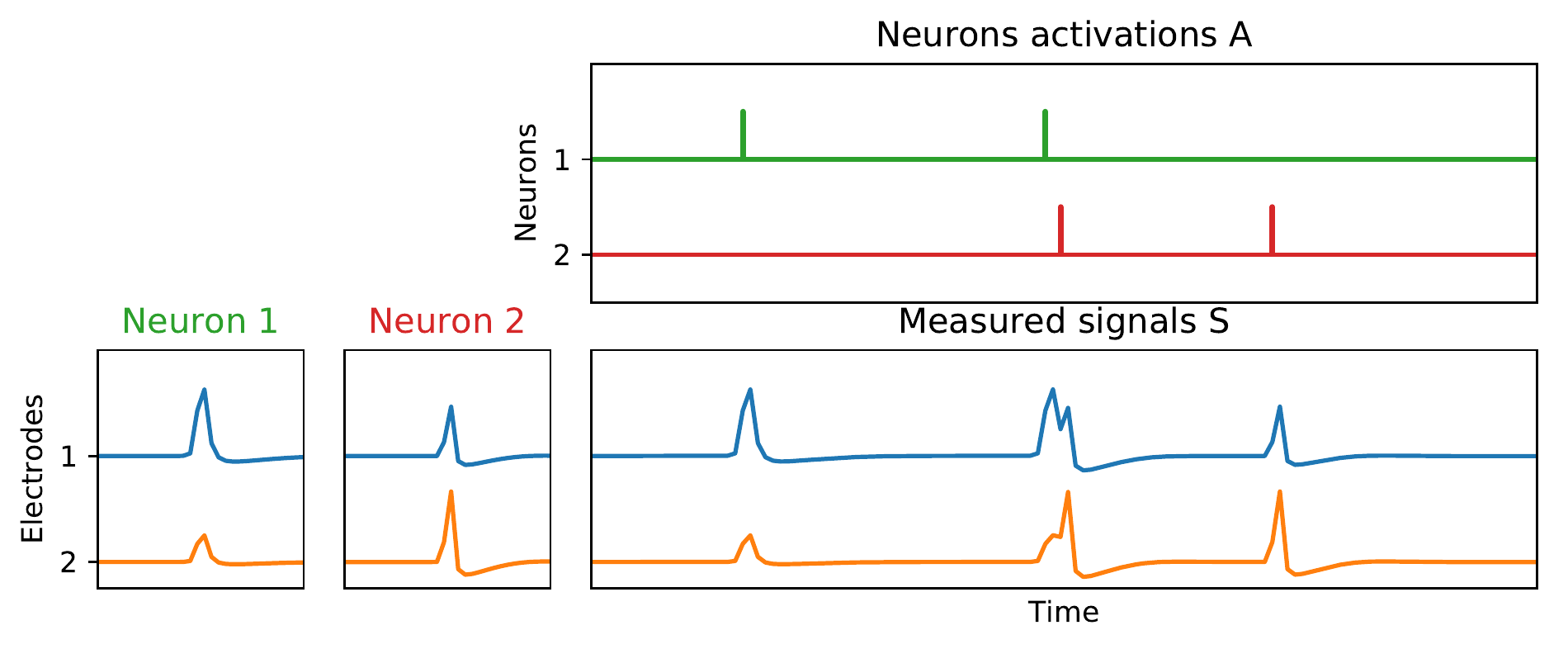}
        \caption{Convolutional model illustration in the simple setting with $E=2$ electrodes and $N=2$ neurons. The shapes of the action potentials are represented in the bottom left part. Remark that for each activation in the top right corner, the corresponding shape appears on the recorded signal in the bottom right corner. Also note the unusual shapes recorded around the center of the graph. This phenomenon may appear when two neurons activate at almost the same time (synchronization). This causes a superposition of the action potentials that creates shapes which are not in the model.}
        \label{fig:convolutional_model}
    \end{center}
\end{figure}

During an experiment, $E$ electrodes record the activity of $N$ neurons. Each electrode records a signal of size $T$ (number of time steps). We propose to model the link between the activity of the neurons and the recorded signals using a convolutional model, introduced by \cite{ekanadham2011recovery}. This model is written as
\begin{equation}
\label{eq:convolutional_model}
	\tens{Y} = \sum_{n=1}^{N} \tens{W}_n * \vec{a}_n^\star + \boldsymbol{\Xi},
\end{equation}
where $\tens{Y} \in \mathbb{R}^{E \times T}$ is the matrix of the observations containing the $E$ recorded signals of size $T$, $\boldsymbol{\Xi} \in \Rn{E \times T}$ is a
random noise matrix, and $*$ is the \textit{convolution} operator along time defined for any two vectors $\vec{x} \in \Rn{\ell}$ and $\vec{y} \in \Rn{T}$ as
\begin{equation*}
    (\vec{x} * \vec{y})_n = \sum\limits_{i=1}^{\min(\ell,n)} x_i \ y_{n-i+1}.
\end{equation*}
Note that we suppose $0$ values for $\vec{y}$ outside of its support as is
classical in signal processing \footnote{For a Python implementation where the handling of the
convolution on the borders can be passed as parameter, see for
instance
\url{https://docs.scipy.org/doc/scipy/reference/generated/scipy.signal.convolve.html}}.
An illustration of the model for given parameters $\tens{W}_n$ and $\vec{a}_n$ is provided in Figure~\ref{fig:convolutional_model}.
The matrix $\tens{W}_n =[\vec{w}_{n,1},\dots,\vec{w}_{n,E}]^\top \in \Rn{E \times
\ell}$ contains the shapes $\vec{w}_{n,e}$ of the action
potentials of neuron $n$ on every electrodes $e$.
Note that each shape is described
by $\ell\ll T$ points. This model assumes that for any neuron/electrode couple,
its shape does not change along time (stationarity of the shape). 
The vector $\vec{a}_n^\star \in \mathbb{R}^{{T}}$ is called the activation vector of neuron $n$.
The non zero entries of $\vec{a}_n^\star$ correspond to the activation times of this
neuron. We recall in table~\ref{table:orders_magnitude} the main quantities from the model and data and provide some orders of magnitude.
Note that while a natural assumption would be  that vector $\vec{a}_n^\star$
is binary (0/1 values), the
resulting optimization problem becomes NP-complete. In addition this does not
allow to model the change of amplitude for the action potentials, which can
occur for neurons \cite{lewicki1998review}, typically after a burst of activity.

Let us write $\tens{A}^\star$ the matrix in \Rn{N\times T} which contains all the $\vec{a}_n^\star$ in its columns. Estimating $\tens{A}^\star$ when the number of neurons is greater that the number of electrodes would be impossible to tackle without additional structural assumptions.
Since the firing rates of the neurons are very small compared to the sampling rate of the signal, $\vec{a}_n^\star$ is clearly a sparse vector. This encourages us to consider an estimator of $\tens{A}^\star$ promoting sparsity. We choose the well-known Lasso estimator, originally proposed by \cite{tibshirani1996regression} and applied to spike sorting problems in \cite{ekanadham2011recovery}. The problem writes as:

\begin{equation}
\label{eq:lasso_convo}
    \min\limits_{\tens{A}}  \quad \norm{\tens{Y} - \sum_{n=1}^{N} \tens{W}_n * \vec{a}_n}{2}^2 + 2 \lambda \sum_{n=1}^N\|\tens{a}_n\|_1,
\end{equation}
where $\lambda>0$ is the regularization parameter of the Lasso and the norm are
defined as the Frobenius norm
$\|S\|_2=\sum_{i,j} S_{ij}^2$ and $\ell_1$ norm $\|\vec{a}\|_1=\sum_{i} |a_{i}|$. It is the only tuning parameter of the method and its
value depends on the
signal-to-noise ratio.

The activations $\vec{a}_n$ and the shapes $\tens{W}_n$ of the model
\eqref{eq:convolutional_model} could be estimated by alternative optimization.
However, due to the orders of magnitude presented in
table~\ref{table:orders_magnitude}, the hardest step would be the update of the
$\vec{a}_n$: indeed, their number of variables growths linearly with $T$, while
$\ell$ is fixed and small. Therefore, as announced in the introduction, we
decide to focus on the estimation of the activations $\vec{a}_n$ while the
shapes $\tens{W}_n$ are assumed to be known. In practice, we can obtain the
shapes and the number of neurons from a classical spike sorting algorithm. Our
estimation of the activations $\vec{a}_n$ then allows to handle correctly the
synchronizations and all the spikes that were not sorted by the first algorithm. \pat{Note also that the biological signal that we are modeling here, is really small so that there is no problem in considering an additive superposition of the shapes since it can never reach the saturation of the device.}

    \begin{table}[t]
        
        \begin{center}        
            \begin{tabular}{|l|c|c|}            
                \hline
                Quantity & Notation & Orders of magnitude \\
                \hline
                Number of electrodes & $E$ & 4-4000 \\
                \hline
                Number of neurons & $N$ & 1-1000 \\
                \hline
                Number of time steps & $T$ & $10^8$ \\
                \hline
                Shape length & $\ell$ & 30-150 \\
                \hline
            \end{tabular}
        \end{center}
        \caption{Orders of magnitude for the quantities of the problem. This table illustrates a typical situation for a recording of one hour at the sampling rate of 30kHz. Taking a realistic neuron spiking rate of 30Hz, we then can expect that the number of non zero coordinates in each $\vec{a}_n$ is about $10^5$, which is very small compared to their size of $10^8$.\label{table:orders_magnitude}}
    \end{table}

\subsection{Vectorizing the Lasso problem and optimality condition}
\label{subsec:vectorization}

\paragraph{Vectorized model}
Since the convolution is a linear operator, the convolutional model
\eqref{eq:convolutional_model} can be translated as a linear model, doing a
vectorization step. We define the vectorization as the concatenation of the
temporal signals, \emph{i.e.} the lines of the multivariate signals, with
$\vec{y}=(Y_{1,1},Y_{1,2},\dots, Y_{E,T-1},Y_{E,T})^\top$ the vectorization of the
recorded signals in $\vec{Y}$ and $\boldsymbol{\xi}=(\Xi_{1,1},\Xi_{1,2},\dots, \Xi_{E,T-1},\Xi_{E,T})^\top$ the
concatenation of the noise signals in $\boldsymbol{\Xi}$ and
$\vec{a}^\star=(A_{1,1}^\star,A_{1,2}^\star,\dots,
A_{N,T-1}^\star,A_{N,T}^\star)^\top=(\vec{a}^{\star \top}_1,\dots,\vec{a}^{\star
\top}_N)^\top$ the concatenation of the activations in $\vec{A}^\star$. Note that in
the following we will sometime index vectors $\vec{a}$ with doubles indices
${a}_{n,t}=A_{n,t}$ for readability reason. The linear
convolutional model in \eqref{eq:convolutional_model} can be expressed in
vectorized format as
\begin{equation}
    \vec{y} = \vec{H}\vec{a}^\star+\boldsymbol{\xi}
    \label{eq:vec_model}
\end{equation}
where the matrix $\vec{H}\in\R^{ET \times NT}$ is a very structured block
Toeplitz matrix (due to the convolutions by the different shapes in $\vec{W}$).
The columns of $\vec{H}$ will be indexed here by a time $t$ and neuron index $n$
for a better readability. The matrix
$\vec{H}=(\vec{h}_{1,1},\vec{h}_{1,2},\dots, \vec{h}_{N,T-1},\vec{h}_{N,T})$
is the concatenation of columns $\vec{h}_{n,t}$ corresponding to an activation
of neuron $n$ at time $t$. The column $\vec{h}_{n,t}$ can be recovered from the
shapes $\vec{w}_{n,e}$ with
$\vec{h}_{n,t}=((\pusht{\vec{w}_{n,1}}{t})^\top,\dots,(\pusht{\vec{w}_{n,E}}{t})^\top)^\top$
where $\pusht{\vec{w}}{t}$ is a vector of size $T$ where the shape $\vec{w}$ has
been pushed at position $t>0$ and its remaining components have value $0$.

\paragraph{Optimization problem and KKT}
The sparse estimation of the activations with the Lasso problem
\eqref{eq:lasso_convo} reformulates after vectorization as 
\begin{equation}
\label{eq:lasso_problem}
    \vec{\hat{a}}
    =
    \argmin_{\vec{a}  \in \Rn{NT}} \quad \norm{\vec{y} - \tens{H}\vec{a}}{2}^2 + 2 \lambda \norm{\vec{a}}{1}.
\end{equation}
An essential property of the Lasso are the following necessary and sufficient
optimality conditions. $\vec{\hat{a}} \in \Rn{NT}$ is a Lasso solution, that is
a solution of problem \eqref{eq:lasso_problem}, if and only if
\begin{equation}
    \label{eq:optimality_conditions}
    \begin{cases}
        \vec{h}_{n,t}^\top (\vec{y} - \tens{H} \vec{\hat{a}}) = \lambda \ \text{sign}(\hat{a}_{n,t}) & \text{, if } \hat{a}_{n,t} \neq 0,\\
        | \vec{h}_{n,t}^\top (\vec{y} - \tens{H} \vec{\hat{a}}) | \leq \lambda & \text{, if } \hat{a}_{n,t} = 0.
    \end{cases}
\end{equation}
In particular for a Lasso solution $\vec{\hat{a}}$, we have for any $1 \leq t
\leq T$ and $1 \leq n \leq N$
\begin{equation}
    \label{eq:kkt_condition}
    \text{ if } \quad	| \vec{h}_{n,t}^\top(\vec{y} - \tens{H} \vec{\hat{a}}) | < \lambda,\quad \text{then } \quad \hat{a}_{n,t} = 0.
\end{equation}

Condition~\eqref{eq:kkt_condition} is simple to test and really useful in a
context of high sparsity where the majority of coordinates in $\vec{\hat{a}}$ are
optimal in $0$. This suggests the use of an iterative scheme for the resolution
of the Lasso problem: starting from the null vector, we can activate iteratively
the coordinates of $\vec{\hat{a}}$. This strategy, called working set, is
presented in more details in the next subsection~\ref{subsec:as}.
    
To describe this algorithm, we need the following notation. For any $J \subset \{1,
\dots, NT\}$ and any $\vec{a} \in \Rn{NT}$, we define $\vec{a}_J$
as the vector obtained by keeping only the coordinates from $\vec{a}$ which are
in $J$. 
We also define $\tens{H}_J$ as the matrix obtained by keeping only the columns
from $\tens{H}$ which are in $J$.

\subsection{Generic working set algorithm}
\label{subsec:as}

From a computational point of view, computing straightforwardly a Lasso
estimator, ie a solution of the problem \eqref{eq:lasso_problem}, can be very
expensive in high dimension even when efficient convolution can be implemented
instead of the full matrix product in \eqref{eq:lasso_problem}. In this work, we
harness the working set strategy (also known as active set) 
\cite{lee2007efficient,szafranski2008hierarchical,boisbunon2014} in order to
compute it more efficiently. The main idea of the
working set is to activate sequentially the coordinates of the solution using the
optimality condition \eqref{eq:kkt_condition} as a criterion.

\paragraph{Principle of the algorithm}
We call working set and we note $J$ the set of the active coordinates in $1,\dots,NT$ of the
solution. We provide a generic formulation of the algorithm in Algorithm
\ref{algo:gas}. Initializing the solution as the null vector, we proceed
iteratively as follows: as long as the optimality condition is not verified, we
activate the coordinate $j_0$ that violates the most the KKT condition \eqref{eq:kkt_condition}
(line~\ref{algo:gas_max} of algorithm~\ref{algo:gas}), we add it to $J$
(line~\ref{algo:gas_add}), and we update the solution by solving the Lasso
problem on the current working set with $\tens{H}_J$ (line~\ref{algo:gas_sublass_solve}). Thanks
to the sparsity of
$\vec{a^*}$, we expect to solve several Lasso problems of size $\leq |J|$, which remains
small compared to $ET$. In the worst case, the working set algorithm would
activate every possible coordinates, thus it ends in finite time.
In practice, the computation of the optimality condition vector
(line~\ref{algo:gas_kkt}) and the update of the solution
(line~\ref{algo:gas_sublass_solve}) are the most expensive steps of the
algorithm. 

\paragraph{Efficient implementation on convolutional models} {One simple
approach to solve problem \eqref{eq:lasso_convo} would be to pre-compute the
matrix $\tens{H}$ and use it for the KKT line \ref{algo:gas_kkt} and Lasso
solvers line \ref{algo:gas_sublass_solve} in Algo. \ref{algo:gas}. But this
approach does not scale properly in memory since the memory complexity of
storing $\tens{H}$ is $O(ENT^2)$ where $T$ is typically very large. The
computational complexity of line \ref{algo:gas_kkt} for a dense matrix
$\tens{H}$ is also $O(ENT^2)$ which does not scale well with the problem
dimensionality. 

 But in practice $\tens{H}$ is very sparse, completely defined through $\tens{W}$, and the convolutional operator can be
 computed exactly with a much smaller memory footprint. Using direct convolution
 instead of a general matrix product, one can compute the gradient of the
 quadratic loss in line \ref{algo:gas_kkt} for a computational complexity of
 $O(ENT\ell)$ and a memory complexity of  $O((E+N)T)$. This means that this
 operation can be used to compute  efficiently the KKT condition in the working
 set and to compute the gradient in the inner lasso solver line
 \ref{algo:gas_sublass_solve}. Also note that in addition to the efficient
 implementation for the convolution operator, the matrix $\tens{H}_J$ 
 in the
 working set is still very sparse with only an order of $O(E|J|\ell)$ non-sparse
 lines, which means that the inner solver can be solved exactly on a much
 smaller subproblem with a matrix of size $O(E|J|\ell)  \times O(|J|)$.

 The  efficient implementation discussed above allows to solve larger problems
 but several computational bottleneck persists: at each iteration in the working
 set one needs to perform  $O(ENT\ell)$ operations and since the number of
 iterations of the working set will be proportional to $T$ it leads to at least a
 quadratic complexity \emph{w.r.t.} $T$, which again does not scale well and does not
 allow to provide real time spike sorting. We propose in section \ref{sec:swas}
 the idea of sliding window working set, which takes advantage of the temporal structure of
 the problem to solve it more efficiently.
 
}

\begin{algorithm}[t]
	\caption{Generic working set algorithm \label{algo:gas}}
	\begin{algorithmic}[1]
		\REQUIRE $\vec{y}, \tens{H}, \lambda>0, \epsilon>0$
		\STATE $J \leftarrow \emptyset$, $\vec{\hat{a}} \leftarrow \mathbf{0}$
		\REPEAT
		\STATE $\vec{g} \leftarrow \tens{H^\top} (\vec{y} - \tens{H} \vec{\hat{a}}) $ \label{algo:gas_kkt}
		\STATE $j_0 \leftarrow \argmax_{l \in J^c} |g_l|$ \label{algo:gas_max}
		\STATE $J \leftarrow J \cup \{j_0\}$ \label{algo:gas_add}
		\STATE $\vec{\hat{a}}_J \leftarrow$ Solve Lasso \eqref{eq:lasso_problem} for
		sub-problem $(\tens{H}_J, \vec{y})$ \label{algo:gas_sublass_solve}
		\UNTIL $g_{j_0} < \lambda + \epsilon$
		\RETURN $\vec{\hat{a}}$, $J$
	\end{algorithmic}
\end{algorithm}

\subsection{Biologically based assumptions}
\label{subsec:bio_assumptions}

We present here some biological properties about neurons and action potentials
and explain how these properties translate into mathematical assumptions related
to our model. Taking advantage of these properties in later sections, we
demonstrate that our estimator of $\vec{a}\kl{^*}$ verifies nice statistical and
computational properties. Moreover it also allows us to derive the theoretical temporal
complexity of our algorithm.

First we present the following assumption about the support  $\Sstar=Supp(\vec{a}^\star)$ of the
true model parameter $\vec{a}^\star$. 
\begin{assumption}\label{ref}
(Absolute refractory period) All indices in the support $\Sstar$ are at least $\l+1$ apart.
\end{assumption}

This assumption is a mathematical reformulation of the idea of refractory period of a neuron. Right after a
neuron fired an action potential, there is a short period of time during which
the neuron can not fire again. In the following, we assume that every neuron in
the model shares the same refractory period. 
As stated in page 4 of \cite{dayan2001theoretical}, this
refractory period length is several times larger than the action potential
length. For simplicity reason we suppose
that this period is  ${\ell}$, which is the length of the potential shape
window. This assumption is of particular importance in our case: it means that the activation of a particular neuron are far away temporally which makes them easier to identify statistically.

After making some assumption on the support of the true model, we  make 
some assumptions on the shapes of the action potentials of the individual
neurons. These shape assumptions are better described as properties of the Gram
matrix $\tens{G}=\tens{H}^\top\tens{H}$. From the definition of
the columns of $\tens{H}$ {(see below Equation \ref{eq:vec_model})}, we can recover the
following Lemma
\begin{lemma}
\label{lem:Gblockband}
     For all $t$ and $t'$ in $\{1,...,T\}$ and for all $n$ and $n'$ in $\{1,...,N\}$, we have:
     \begin{equation}
        G_{(n,t),(n',t')}=\vec{h}_{n,t}^\top\vec{h}_{n',t'}=\sum_{e=1}^E   
        (\pusht{\vec{w}_{n,e}}{t})^\top\pusht{\vec{w}_{n',e}}{t'}
                 \label{eq:gram_comp}
     \end{equation}
    In particular it is null when $|t-t'|>\l$.
\end{lemma}

Note from the Lemma above that similarly to the columns
$\vec{h}_{n,t}$ that are indexed by neuron $n$ and time $t$ we will index the
components of $\tens{G}$ such as
$G_{(n,t),(n',t')}$. 
We now define below several correlation
assumptions on the shapes.

\begin{assumption} \label{onG} ~\\
{\bf 2.a} (Neurons are recognizable) There exists $\epsilon>0$ such that for all
$n\neq n'$ and $|t-t'|\leq \l$
\begin{equation}
    |G_{(n,t),(n',t')}|\leq \epsilon.
    \tag{2.a}
    \label{eq:2a_neurons}
\end{equation}
{\bf 2.b} (Spikes are peaky) There exists $\rho \in (0,1)$, such that for all  $n$ and
$0<|t-t'|\leq \l$
\begin{equation}
    |G_{(n,t),(n,t')}|\leq \ld{\rho |G_{(n,t),(n,t)}|}.
    \tag{2.b}
    \label{eq:2b_peaks}
\end{equation}
{\bf 2.c} (Shapes have bounded energy). There exist $\overline{c}>\underline{c}>0$ such that
for all  $n$ and $t$,
\begin{equation}
    \underline{c}\leq |G_{(n,t),(n,t)}|\leq \overline{c}.
    \tag{2.c}
    \label{eq:2c_bounded}
\end{equation}
\end{assumption}

Assumption {\bf 2.a} is of great importance for the statistical analysis of our
methodology. It essentially means that the shapes of two distinct neurons are
distinguishable, allowing to attribute each spike to the correct neuron.
  This is
reasonable due to the difference between neurons but also due to their spatial
localization that will means different impact on different electrodes {\cite{biffi}}.
Assumption {\bf 2.b} is also reasonable due to the fact that action potentials
are also called spikes that will definitely diminish their autocorrelation in
the presence of a temporal delay {(see classical shapes for instance in \cite{pouzathuxley})}. Finally Assumption {\bf 2.c} is also
physically plausible since an action potential with too small energy would be
indistinguishable from recording noise and the potential obviously has a bounded
energy \cite{pouzathuxley}. \ld{In order to illustrate the statistical feasibility of our approach, we provide in section \ref{sec:num_exp} numerical values for the parameters $\epsilon, \rho, \underline{c}$ and $\overline{c}$ which we have introduced in these assumptions.Note that these values highly depend on the recordings and should be understood as just possible values.}


\section{Sliding window working set algorithm}
\label{sec:swas}

In this section we present our novel working set algorithm. This algorithm
builds on the fact that thanks to the structure of the problem and of its
solution, it can be decomposed into smaller problems. This will be illustrated
and discussed next in \ref{subsec:overlaps}.  The algorithm is then 
introduced and illustrated in subsection \ref{subsec:presentation_swas} and
discussed more in detail in subsection \ref{subsec:algorithm_math}.

{In the sequel we need the following notation.}
We define a temporal window $\omega = \llbracket \omega_1, \omega_2
\rrbracket$ where $\omega_1 \leq \omega_2$ and $\omega_1,\omega_2\in \{1,\dots,
T\}$ that contain all samples whose temporal index $\omega_1\leq t\leq \omega_2$.
This temporal window will be used in the following to index vectors with 
$\vec{a}_\omega$, that contains the temporal samples  $\omega_1\leq t\leq
\omega_2$ for all neurons $n$, and matrix $\tens{H}_\omega$ where are selected
only the columns $\vec{h}_{n,t}$ where $\omega_1\leq t\leq
\omega_2$, for all neurons $n$.

\subsection{Overlaps and independent problems}
\label{subsec:overlaps}

We introduce here the notions of spatial and temporal overlaps that will be
useful in the remaining. These overlaps will allow us to split the large
optimization problem \ref{eq:lasso_convo} into several smaller scale problems
that are individually easier to solve. We first discuss the notion of spatial
overlap that will be important in the MEA case and is related to the physical
position of the neurons. Next we discuss {the temporal overlaps} that are related to
the temporal activations of the individual neurons.

\subsubsection{Spatial overlaps}
\label{subsubsec:spatial_overlaps}

{
In the MEA case, solving the problem on {the full set of $N$ neurons}  has a heavy
computational cost. In this section, we want to take advantage of an important property
of the problem: the spatial distribution of the neurons. Simply put, we harness
the fact that two physically distant neurons are not recorded by the same
electrodes. Thus 
their respective spikes should not overlap on any electrode, even if these
neurons emit simultaneous spikes. 

{The problem is in fact a bit more complex than that because there might be transitive effect. Indeed in Figure \ref{fig:overlaps}, neuron $N_1$ and $N_3$ are recorded by disjoint sets of electrodes, but still it is not possible to speak of disjoint independent Lasso problems, because their spikes might be mixed with the ones of $N_2$. Hence we need to access the spatial overlaps that they form. We hope that these overlaps will form much smaller sets that the complete set of neurons and this is linked to the spatial distribution of the neurons. So let us first precise why such 
a phenomenon might appear from a biological point of view in the MEA case. In the tetrode case, the number of electrodes is so small that such phenomenon is not relevant.}

\paragraph{Neuron density and localization}

Typical studies of {\it in vitro} cultures report roughly 1000  neurons per
$mm^2$ \cite{biffi}. Note that these cultures usually provide a higher density
of cells than in {\it  ex vivo} experiments where slices of brain are used. The
range between electrodes in a MEA depends on the type of MEA and might range
from 200 $\mu m$ \cite{biffi} to about $30 \mu m$ \cite{muthmann}.  
Finally the electrical signal that is emitted by a neuron suffers from various
kind of attenuation and people analyzing MEA signals usually think that a neuron
is recorded by very few nearby electrodes (for instance only 5 electrodes in the MEA are used by
\cite{muthmann}, which corresponds to a range of about  200 $\mu m$). These orders of magnitude mean
that in practice the impact of the activation for a given neuron will be very
localized between a few electrodes, which introduces nice properties discussed below.
{More precisely, in Section \ref{subsubsec:results_spatial_overlaps}, we will
leverage percolation results to upper bound  the size of the spatial overlaps
with large probability.}

To fix the configuration, from now on, the MEA case corresponds to $E$ electrodes placed on a square lattice.

\paragraph{Spatial clustering of the neurons}
Let us now formalize the concept of spatial overlaps to explain the algorithm.
Two neurons $n$ an $n'$ are independent when we have:
\begin{equation}
	 \forall e, \quad \vec{w}_{n,e}=0 \ \text{or} \ \vec{w}_{n',e}=0.
	\label{eq:indep_neurons}
\end{equation}
This condition is true when one of the two shapes
is the null vector which happens when the two neurons do not share any electrode
where they are both recorded (they do not overlap). Using this pairwise
independence, one can easily construct a clustering of the neurons as illustrated in Figure
\ref{fig:overlaps}.left where 3 independent spatial overlaps are recovered. Note that
\eqref{eq:indep_neurons} implies that between two neurons $n$ and $n'$ in two
independent clusters we have $\vec{h}_{n,t}^\top\vec{h}_{n',t'}=0,\quad \forall
t,t'$, which means that both the quadratic term and the Lasso regularizer are
separable in several Lasso subproblems {(one by spatial overlap) and that they} can be solved
independently. This means that by performing beforehand a clustering of the neurons based on the shapes of their action potentials, we can
greatly decrease the complexity of the problem. In the following we will suppose
that this clustering has been done and that the problem is solved on a subset of
neurons (in a spatial overlap) and of electrodes (the ones active inside the spatial overlap).  }

\subsubsection{Temporal overlaps}
\label{subsubsec:temporal_overlaps}

\begin{figure}
	\begin{center}
		\includegraphics[width=340pt]{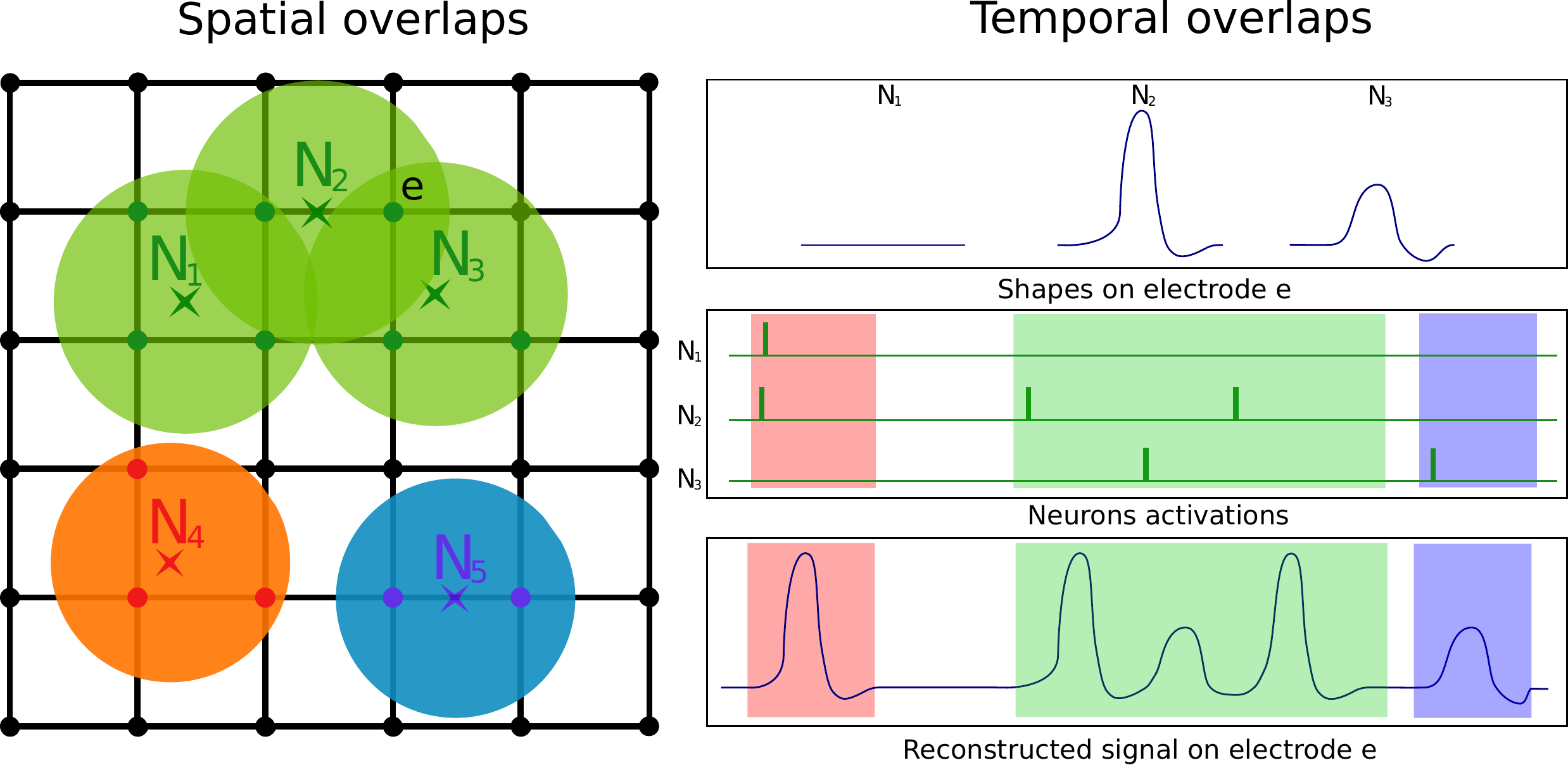}
		\caption{Illustration of the spatial and temporal overlaps. At the left
		hand side, we present an example of 3 spatial overlaps in the case of 5
		neurons, on a regular grid of 36 electrodes. The position $N_j$ of
		neuron $j$ is represented by a check, and the reach of its spikes by a
		disc of radius $r$.  On the right
		hand side, we provide an example of temporal overlaps for the neurons 1,
		2 and 3. We provide the shapes of each neuron and the reconstructed signal on the
		electrode $e$. Remark that since neuron 1 is far away from $e$, its
		shape on $e$ remains at 0. The independent spatial and temporal overlaps
		are illustrated with different colors.
		\label{fig:overlap} \label{fig:overlaps}}
	\end{center}
\end{figure}

{
In addition to splitting the optimization problem thanks to the spatial overlaps
of the neurons, one can also use the structure of the problem to split the
problem into almost independent temporal windows. Let us first review the biological phenomenon which explains such reduction.

\paragraph{Neuron activations and refractory period}
Neurons fires quite scarcely and some classical models are either
 Poisson processes in continuous time with a frequency of usually 10Hz (max 100Hz)
 or their discrete counterpart that are Bernoulli process
(see for instance \cite{tuleau,goodness} and the reference therein).
 Both models have been adapted to encode the refractory period, for instance
 using Poisson with
 dead time, which basically consists in erasing the spikes that are
 too close. These more precise variations can only decrease the number of
 spikes. In case of synchronization, the synchronizations between
 neurons are usually simulated by joint Poisson or Bernoulli process  \cite{tuleau}, so that
 the firing pattern of the whole system remains globally Poisson (or Bernoulli).

\paragraph{Temporal overlap and independent windows} Similarly to spatial
overlaps we can find independent temporal clusters (temporal windows) of
activations. We define two activation at times $t$ and $t'$ as independent if
$|t-t'|>\ell$. Indeed in this case the supports of the convolution (of size $\ell$)
do not overlap and it is easy to show that
$\vec{h}_{n,t}^\top\vec{h}_{n',t'}=0,\quad \forall n,n'$. This is interesting
because it means that for any windows $\omega$ and $\omega'$ such that $\omega_2
< \omega_1'+\ell$ we have $\tens{H}_\omega ^\top\tens{H}_{\omega'}=\mathbf{0}$
where $\mathbf{0}$ is the null matrix. This implies again that the optimization
problem can be solved independently on $\omega$ and $\omega'$.

Similarly to spatial overlap, one can find independent windows that contain the
activations of the neurons, as illustrated in Figure \ref{fig:overlaps}.right.
But note that this time the temporal overlaps
cannot be found \emph{A priori} since the
actual support of the temporal activations is not known. This means that despite
this nice separability of the problem, one cannot use it to speedup the
optimization until the support of the solution is known. The main motivation for the
sliding window working set algorithm introduced below is to find this support
and independent windows in an efficient and online way.
}
{Mathematically, the size of the temporal overlaps themselves is also controlled with large probability (see Section \ref{subsubsec:results_spatial_overlaps}) and this will impact the overall complexity of the algorithm.}

\subsection{Sliding window working set for the global problem}
\label{subsec:presentation_swas}

{We present here our sliding window working set algorithm. The main idea of
the algorithm is to work only on a small temporal window and use the working set
principle to simultaneously solve the optimization problem in the window and
find the window that is independent from the rest of the signal.
Note that this
algorithm is used on each independent spatial overlap
so in
fact for small values of $N$ and $E$ in the MEA case. 

\paragraph{Principle of the algorithm} The main algorithm is detailed in
Algorithm \ref{algo:swas_detailed} where line 3 denotes an update of the
large vector $\hat{\vec{a}}$ on the current window. The idea is to solve the Lasso problem on
small windows $\omega$ starting with the beginning of the signal $\omega
= \llbracket 1, 4\ell\rrbracket$ 
and perform the following operations until the
end of the signal is reached:
\begin{enumerate}
	\item Computing the Lasso solution $\vec{\hat{a}}_\omega$, on the  window
	$\omega$ with the working set algorithm.
	\item Updating the window $\omega$ depending on the support of $\vec{\hat{a}}_\omega$:
	\begin{enumerate}
		\item If the support $Supp(\vec{\hat{a}}_\omega)\in \llbracket \omega_1+\ell,
		\omega_2-2\ell \rrbracket$ then the current problem is independent from
		the rest of the signal and the window is updated as $\omega
		= \llbracket \omega_2+1-\ell,\omega_2+3\ell		\rrbracket$.
		\item If the support $Supp(\vec{\hat{a}}_\omega) \cap \llbracket\omega_1,
			\omega_1+\ell-1 \rrbracket \neq \emptyset$ has components in the first $\ell$ samples of the window, then we merge the current window with the last (because the KKT conditions on the last $\ell$ samples in the previous window have changed.)
		\item Else the window needs to be extended as $\omega = \llbracket \omega_1, \omega_2+ \ell \rrbracket$
	\end{enumerate}
		\end{enumerate}

Once the Lasso is solved on the window $\omega$  in step 1, the optimality
conditions are verified on the window. If the support of the activation $Supp(\vec{\hat{a}}_\omega)\subset \llbracket \omega_1+\ell,
\omega_2-2\ell \rrbracket$, it means that the reconstructed signal after
convolution  signal is entirely contained into $\omega$. 
This means that the solution on this window is independent from the previous one
and probably independent from the next one (since there is no activations in the
last $\ell$ samples of the window).  In this case, we have found the Lasso
solution for the previous window, and then we can work on a new window
immediately after it (step 2.a). \rf{The right border was chosen as $2\ell$
instead of $\ell$ to promote a better exploration ahead of time and well separated temporal clusters
minimizing the occurrence of the more expensive case 2.b discussed below.}
If there are activations at the beginning of the window however (first $\ell$ samples), it means that the residual and the KKT conditions have changed on the last $\ell$ samples of the previous window and we need to potentially update the model there, so we merge the current and previous windows. 
Else, if there are activations in the last $\ell$ samples of the window, we
extend it in order to ensure that at least $\ell$ samples are not activated at the end. Finally for each iteration after updating the current window, we solve again the Lasso on this window efficiently thanks to the working set strategy.  For illustration, one execution of the algorithm
with one electrode and two neurons is provided in Figure
\ref{fig:algo_illustration}. It shows both configurations: when the window is
extended and when the window is shifted to the right.

\begin{algorithm}[t]
	\label{algo:swas_detailed}
	\caption{Sliding window working set}
	\begin{algorithmic}[1]
		\REQUIRE $\vec{y}, \tens{H}, \lambda>0$
		\STATE $\vec{\hat{a}} \leftarrow \mathbf{0}$, 
		 $\omega =  \llbracket \omega_1, \omega_2 \rrbracket
		\leftarrow \llbracket 1, 4\ell \rrbracket$, Empty list of windows $\Omega=[]$
		\REPEAT
		\STATE $\vec{\hat{a}}_\omega \leftarrow$ Solve Lasso with algo. \ref{algo:gas} for sub-problem $(\tens{H}_{{\omega}} ,
		\vec{y})$ using warm-start $\vec{\hat{a}}_\omega $ \label{algo:swas_sublass_solve}
			\IF{$Supp(\vec{\hat{a}}_\omega)\subset\llbracket \omega_1+\ell,
			\omega_2-2\ell \rrbracket$} \label{algo:swas_if}
			    \STATE $\Omega\leftarrow [\Omega, \omega]$ \hfill//Insert current window $\omega$ at the end of list $\Omega$
				\STATE $\omega
		\leftarrow \llbracket \omega_2+1-\ell,\omega_2+3\ell \rrbracket$
		\hfill//Independent 
problem solved so move window
		to next time segment 
		\ELSIF{$Supp(\vec{\hat{a}}_\omega) \cap \llbracket\omega_1, \omega_1+\ell-1 \rrbracket \neq \emptyset$}
		 \STATE {$\tilde \omega , \Omega \leftarrow$ Return last window $\tilde\omega=\llbracket \tilde\omega_1,\tilde \omega_2\rrbracket$ in $\Omega$ and remove it from the list $\Omega$}.
		 \STATE $\omega
		\leftarrow \llbracket \tilde\omega_1,\omega_2\rrbracket$ 	\hfill// merge current window with last window
						\ELSE
			\STATE $\omega
		\leftarrow \llbracket \omega_1,\omega_2+\ell \rrbracket$ \hfill // Extend
window to find the independent temporal overlap		\label{algo:swas_fenetre_max}
			\ENDIF 		\UNTIL $\omega^{(m)}_2 \geq T$
		\RETURN $\vec{\hat{a}}$, $\Omega=[\omega_1,\omega_2,\dots]$
	\end{algorithmic}

																				\end{algorithm}

\paragraph{Algorithm solution \emph{w.r.t.} the original Lasso}
\label{subsec:algorithm_math}
We now address the following question: is the proposed algorithm actually solving the
global optimization problem \eqref{eq:lasso_convo}? We provide to this end the following theorem.

\begin{theorem}
\label{thm:swas_solution_is_lasso_solution}
	The solution $\vec{\hat{a}}$ computed by the sliding window working set is a
	solution of the initial Lasso problem \eqref{eq:lasso_problem}. 
\end{theorem}

\begin{proof}
	By construction of the algorithm, line \ref{algo:swas_if} of the algorithm will return a list of windows $\Omega$ such that $\forall \omega \in \Omega$  the support $Supp(\vec{\hat{a}}_\omega)\in \llbracket \omega_1+\ell,
	\omega_2-2\ell \rrbracket$ which means that the current model will have an
	effect only inside $\omega$ and that for two consecutive windows   $\omega$ and  $\omega'$ in $\Omega$ the temporal indexes  of all active
	variables  $\hat a_t \neq 0$ with $t\in\omega $ and  $\hat a_{t'} \neq 0$ with 
	$t'\in\omega' $ are by construction $\rf{|t'-t|>3\ell}$. This means that as discussed
	in subsection \ref{subsubsec:temporal_overlaps}, the Lasso problems can be
	split as two independent problems on the support and all the other components that are not active respect the KKT meaning that they will be $0$. The solution  $\vec{\hat{a}}$ is
	obtained by successive juxtaposition of the 
	solutions of independent problems on disjoint windows that are estimated
	line \ref{algo:swas_sublass_solve}. 
\end{proof}

\paragraph{Numerical complexity and efficient implementation} We now discuss why
the proposed algorithm is more efficient than the convolutional working set
discussed in section \ref{subsec:as}.
In the standard working set, we recall that the computational cost for the
optimality condition is
of order $O(ENT\ell)$ at each step because of the multiple convolutions
necessary to compute the KKT conditions. In the sliding window working set, the
KKT are computed only on the window $\omega$, and their complexity is reduced to $O(|\omega|NE\ell)$ with
$|\omega| << T$. Proving
mathematically this reduced computational complexity is one of the main focus of
section~\ref{sec:math_results}. Note that the complexities above are given
on the whole problem, but as discussed above, in the MEA case, the spatial clustering of the
neurons means the the complexity depends on the size $E_c$ and $N_c$ of the spatial overlap
 instead of $E$ and $N$.

\paragraph{Comparison with a similar approach} We can find in \cite{moreau2018dicod} a similar strategy, named Distributed Convolution Coordinate Descent (DICOD), for the efficient resolution of large scale convolutional sparse coding formulated as Lasso problems.
By considering a temporal partition of the whole signal, their approach aims at solving small local problems using coordinate descent algorithm. Since the update of a coordinate only influences its vicinity, these local problems can be treated in parallel, almost as independent problems.
In the same manner as theorem~\ref{thm:swas_solution_is_lasso_solution}, they could prove under mild assumptions that the computed solution is indeed a Lasso solution.
Although they demonstrated an important speedup with respect to the global coordinate descent, they do not provide the theoretical complexity of their algorithm with respect to the sizes of the problem.
Taking advantage of the biological assumptions presented in section~\ref{subsec:bio_assumptions}, we prove in the following sections not only that the Lasso estimator retrieves the true support, but also ascertain the theorical complexity of the sliding window working set algorithm. 
Moreover, in our approach, the temporal exploration of the signal with a window allows to treat the problem in an online manner.

\begin{figure}
    \centering
    \includegraphics[width=1\linewidth]{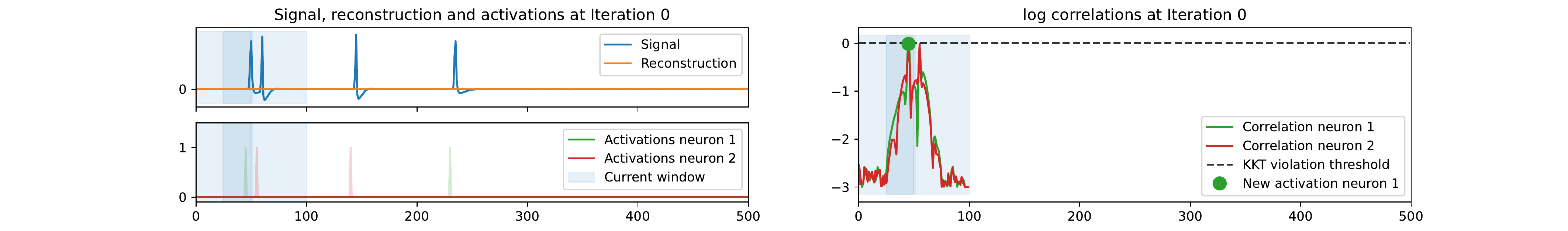}
    \includegraphics[width=1\linewidth]{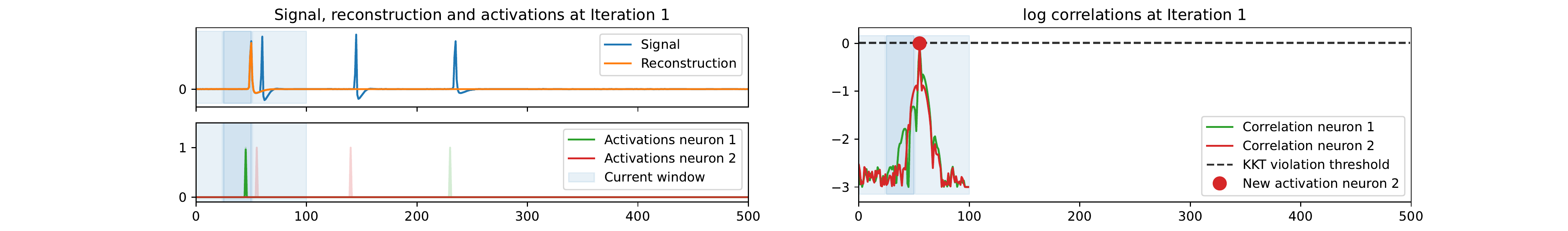}
    \includegraphics[width=1\linewidth]{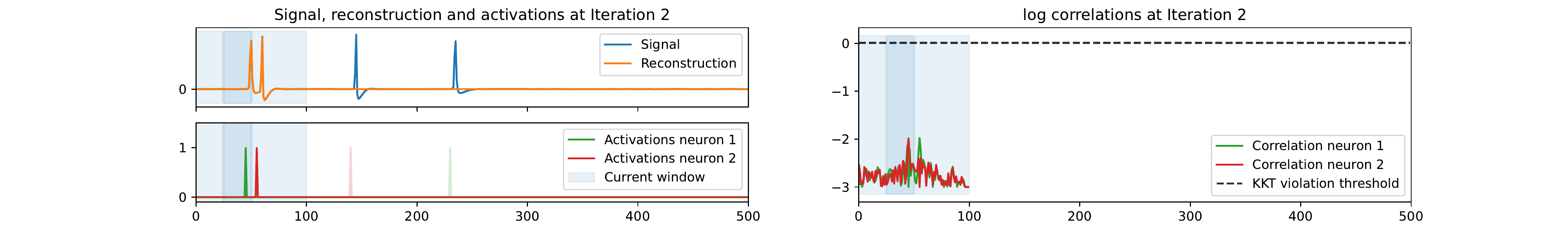}
    \includegraphics[width=1\linewidth]{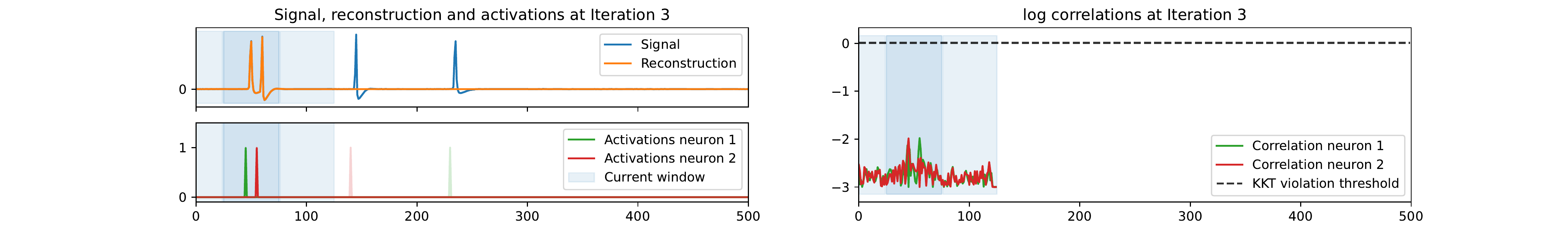}
    \includegraphics[width=1\linewidth]{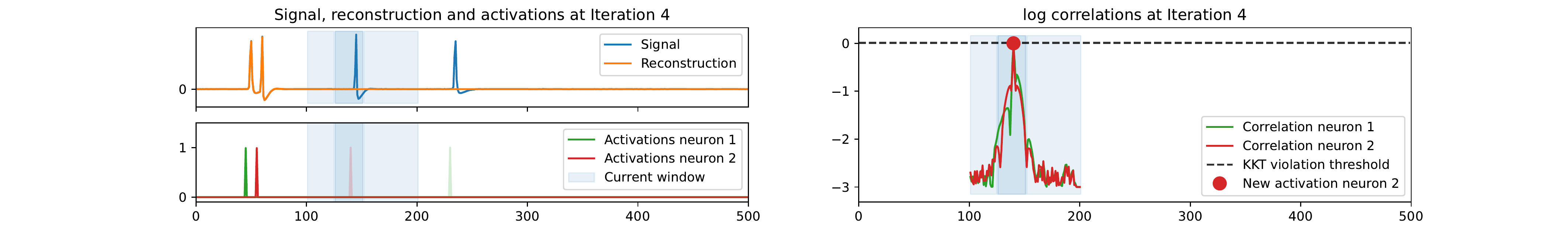}
    \includegraphics[width=1\linewidth]{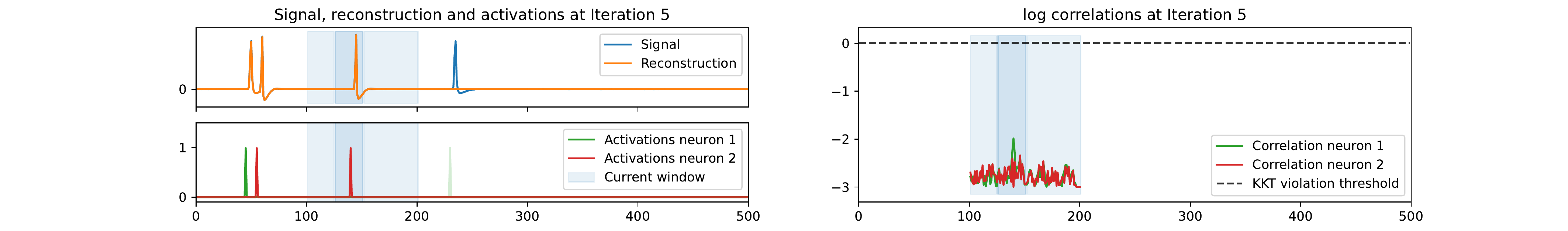}
    \includegraphics[width=1\linewidth]{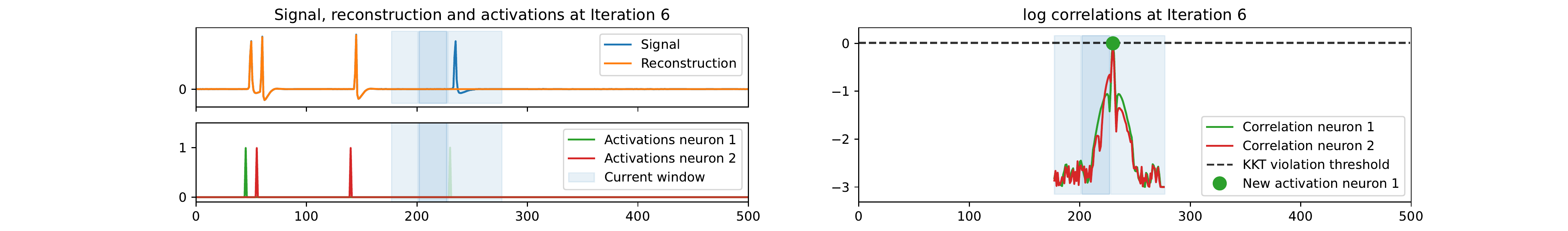}
    \includegraphics[width=1\linewidth]{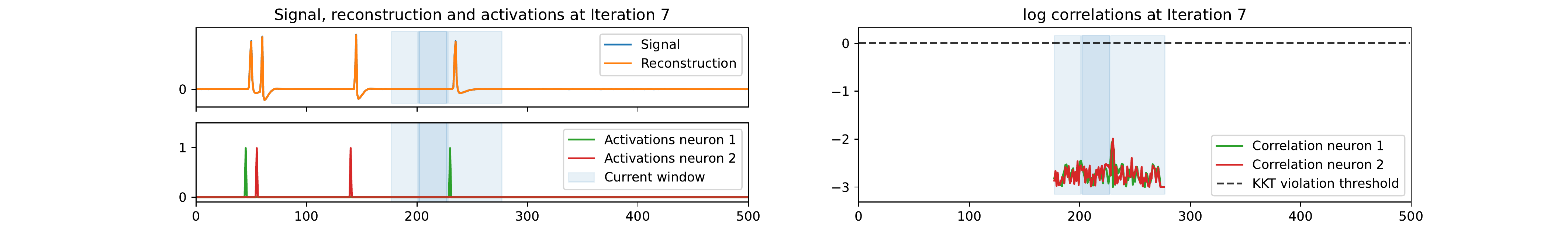}
    \includegraphics[width=1\linewidth]{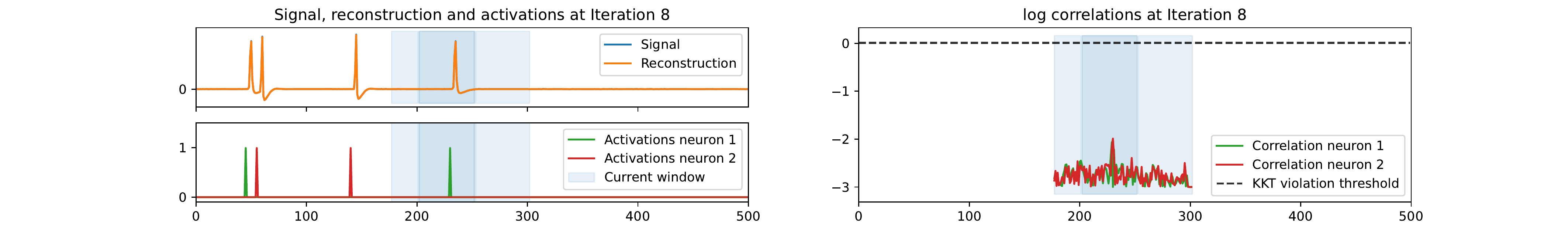}
    \includegraphics[width=1\linewidth]{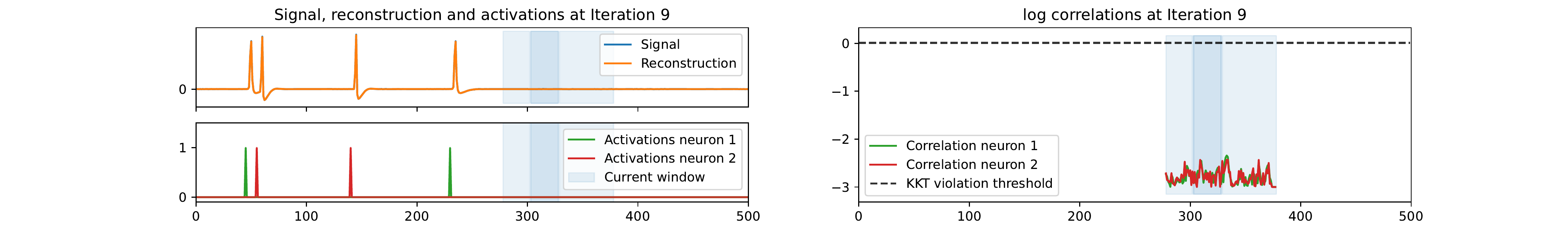}
    \caption{Illustration of the different steps in the proposed algorithm. (left) Observed signal $\vec{S}$ with model reconstruction and sparse model $\vec{a}_i$ The current window is illustrated as a light blue background. True activation are illustrated with transparency (right) KKT violation at the current step. temporal instant and neuron violating the KKT are over the black line. }
    \label{fig:algo_illustration}
\end{figure}

}

\section{Mathematical results}
\label{sec:math_results}

\subsection{Control of the spatial and temporal overlaps}

 \subsubsection{Spatial overlaps \label{subsubsec:results_spatial_overlaps}}
 
 In the MEA case, the $E$ electrodes are placed on a square lattice \rf{with fixed distance $\delta$ between adjacent electrodes}. The range of detection of  a neuron by an electrode is $r_0$. We classically approximate the spatial distribution of the neurons on the lattice by a Poisson process of constant intensity $\gamma$.  
 An electrode detects a neuron if it is at distance less than $r_0$. Therefore two neurons can be detected by the same electrode if their distance is less than $r_0$. If this is the case, we say that these neurons are "connected".  Spatial clusters are given by  maximal sets of neurons that are connected together, or for which there exists a path in between of "connected  neurons".
 
 This framework is known in probability as a particular case of the Poisson-Boolean percolation (see \cite{duminilcopin2018subcritical} and references therein). Thanks to this, we can prove the following proposition.
 
 \begin{proposition}
 \label{spatial}
  There exists a critical value $\gamma_c>0$ which only depends on  $r_0$, such that,  if $\gamma<\gamma_c$,
 then for all $\alpha$ in $(0,1)$, such that $E\geq \tau  \log(3/\alpha)$ for some positive constant $\tau$, there exists an event $\Omega_{\alpha,s}$ of probability larger than  $1-\alpha$, such that on $\Omega_{\alpha,s}$, any spatial overlap  of neurons $c$, with cardinality $N_c$, satisfies
 $$N_c \leq \kappa [\log(E/\alpha)]^2,$$
 with $\kappa>0$, which only depends on $\gamma,\delta, \tau$ and $r_0$.
 \rf{Similarly the number of active electrodes for a given spatial overlap $c$ can be bounded as $E_c\leq \bar\kappa [\log(E/\alpha)]^2$, with $\bar{\kappa}>0$, which only depends on $\gamma,\delta, \tau$ and $r_0$.}
 
 The event $\Omega_{\alpha,s}$ only depends on the position of the neurons on the lattice representing the MEA. 

 \end{proposition}
The proof is given in Section \ref{proof_spatial}.

Let us comment this result qualitatively. First of all $\gamma_c$ is a critical parameter of the percolation theory. When the parameter $\gamma$ is small with respect to $\gamma_c$, as usual for crital percolation parameters, clusters cannot reach infinity, whereas they can if $\gamma$ is too big. As far as we know, precise knowledge of $\gamma_c$ is unknown, but one can still have the following heuristic reasoning : if a neuron can be detected at a range of $r_0$ and if the intensity $\gamma$ (that is, informally, the density of neurons) is very low, it will be quite rare to have two connected neurons, and the cluster size will be roughly one. On the other hand if $\gamma$ is too large, the distance between neighbors will be very small and eventually all neurons will belong to one giant cluster. 

Note that when the shape of the action potential as {perceived} by the various electrodes are known, it is very easy to find these {spatial overlaps} before hand and we will easily know if we are in a subcritical regime where the size of the {spatial overlaps} vary logarithmically  with $E$ or not.

In the rest of the paper, 
\begin{itemize}
\item[(i)]either we focus on a tetrode like case where $E$ is small, so that we discard totally the dependence in $E$, 
\item[(ii)] or on a supercritical regime for a lattice MEA, and then $N_c$ is roughly of the size of $N$, that is the total number of neurons and we can as well solve the big system,
\item[(iii)] or we are in a subcritical regime for a lattice MEA and once restricted to the event $\Omega_{\alpha,s}$, we can solve independently the Lasso problems for each of the spatial overlaps $c$. In this case, the number of neurons is roughly of the order $(\log E)^2$.
\end{itemize}

Note that in (i) or (ii), $N$ is thought to be fixed and a parameter of the problem whereas in (iii) the number of neurons is a random variable and the event $\Omega_{\alpha,s}$ to which we restrict ourselves depends on it.

\subsubsection{Temporal overlaps \label{subsubsec:results_temporal_overlaps}}
{Here we assume that the activations $\vec{a}^*$ are the realisation of a given random process. More specifically, and as explained in Section \ref{subsubsec:temporal_overlaps}, we do not need to model each neuron individually and we do not need to model the amplitude of $\vec{a}^*$. Hence the following formalism can be seen as a worst case scenario.}

We denote $A$ the joint process of length $T$ with values $0$ or $1$, $1$ meaning that at least one of the recorded neurons has fired. Note if we are in the  subcritical regime, $A$ is restricted to the neurons in a given spatial overlap.

We model $A$ by a Bernoulli process of rate $p=N m\Delta$, with $m$ the average firing rate and $N$ the number of neurons (of the spatial overlap possibly) (that is the $A_i$'s are i.i.d. Bernoulli with parameter $p$). Note that in this set up, we force $\Delta<<1/N$ so that we cannot analyze too much neurons at the same time. Another way to see this  is to say that $p$ should be small and to fix ideas we assume that $p\leq 1/2$.

We are saying that two successive spikes $t$ and $t'$ in $A$ are overlapping if $|t-t'|\leq \eta= 4\ell$. 
As from a spatial point of view, from a temporal point of view, the spikes in $A$ that includes all the activation times of all the neurons (of the spatial overlap) are therefore partitioned in overlaps. We can, as for the spatial overlaps, control their size.

\begin{proposition}
\label{tempo}
The temporal overlaps are controlled as follows.
\begin{itemize}
\item In the non-subcritical MEA case or in the tetrode case, with a global activation rate $p=Nm\Delta\leq 1/2$, 
there exists an event $\Omega_{\alpha,A}$ of probability larger than $1-\alpha$ such that on $\Omega_{\alpha,A}$, each temporal overlaps has a length bounded by
$$\mathscr{W}=c' \log(T/\alpha),$$
with $c'>0$ depending only on $\eta$.

\item In the subcritical MEA case, with an activation rate per cluster $p_c=N_c m \Delta\leq 1/2$, 
there exists an  $\Omega_{\alpha,A,s}$ such that, on $\Omega_{\alpha,A,s}$, for each spatial overlap and for each temporal overlap inside a spatial overlap, the size of this temporal  overlap is bounded by
$$\mathscr{W}=c'' \log(E T/\alpha),$$
with $c">0$ depending only on $\eta, \delta, \gamma, \tau$ and $r_0$.
\end{itemize}

\end{proposition}
The proof is given in Section \ref{proof_tempo}.

In the subcritical MEA Case, the condition is on the activation rate per cluster ($p_c\leq 1/2$), which means that the result holds valid even when the global activation rate is $p>1/2$. {In this sense, if the problem is subcritical, even if the MEA is very large and record a huge number of neurons, the size of the  temporal overlaps, which governs the numerical complexity, will still be reasonable.}

Note that in the first case, the event  $\Omega_{\alpha,A}$ depends only on the distribution of the spikes, whereas in the second case, $\Omega_{\alpha,A,s}$ depends both on the spiking distribution but also on the spatial distribution of the neurons.

 \subsection{Control of the noise}

\noindent

\begin{lemma}\label{conc}
Assume that the noises ($(\xi_{e,t})_{e,t}$) are i.i.d. normal random variables with zero mean and finite variance $\sigma^2$. For $\alpha\in (0,1)$, define 
$$
\level_{\alpha} = \sqrt{2\sigma^2 \overline{c} \log\left( \frac{2 N T}{\alpha}  \right)},
$$ 
where $\co$ is given in (\ref{eq:2c_bounded}) and the event
\begin{align}\label{eq:noiseevent}
    \Omega_{\alpha,\noise} = \bigcap_{n,t} \left\lbrace \left| h_{n,t}^\top \noise \right| \leq \level_{\alpha} \right\rbrace. 
\end{align}
 Then we have 
$$
\mathbb{P}\left( \Omega_{\alpha,\noise}\right) \geq 1-\alpha.
$$
\end{lemma}
{The proof is given in Section \ref{proof:lemmaconc}. 

This lemma is very classical and help us to control the level of the noise.
From now on, $\Omega_\alpha$ refers to the event of probability controlled by $1-\alpha$ which is either $\Omega_{\alpha/2,A}\cap \Omega_{\alpha/2,\noise}$ in the tetrode or non-subcritical MEA case, or $\Omega_{\alpha/2,A,s}\cap \Omega_{\alpha/2,\noise}$ in the subcritical MEA case.
}

\subsection{Theoretical properties of the Lasso estimator}
\label{subsec:theoretical_properties_lasso}

\begin{theorem}\label{supp}
Fix $\alpha\in (0,1/2)$.
Let Assumptions \ref{ref} and \ref{onG} be satisfied and let us assume that the noises are i.i.d. Gaussian. 
With the notation of Propositions \ref{spatial}, \ref{tempo} and Lemma \ref{conc},  there exists an event $\Omega_\alpha$ of probability larger than $1-\alpha$ such  that, on $\Omega_\alpha$,  for all temporal window $\omega$ and any solution $\hat{\vec{a}}_\omega$ of the Lasso problem $(\tens{H}_\omega,\vec{y})$ with regularization parameter $\lambda$ (possibly restricted in the subcritical MEA case to any spatial overlap), the following holds.

\begin{enumerate}
\item No spurious activation is discovered, that is
$$
Supp(\hat{\vec{a}}_\omega)\subset \Sstar\cap \omega,
$$
where $\Sstar = Supp(\astar)$ is the true set of activations, as long as
\begin{equation}\label{minlambda}
\lambda > \frac{\cd+ 2\rho \kl{\overline{c}} }{\cd-2\rho\kl{\overline{c}}  - 4\epsilon \mathcal{N}} \left( z_{\alpha} +  2 (\rho \kl{\overline{c}}  + \epsilon \mathcal{N}) \|\vec{a}^*\|_{\infty,\partial \omega}\right) \quad \mbox{and} \quad  \cd>2\rho \kl{\overline{c}} +4\epsilon \mathcal{N}
\end{equation}
with the convention that 
$$\|\vec{a}^*\|_{\infty,\partial \omega}=\sup_{n, t\in \partial\omega} |\vec{a}^*_{n,t}|,$$
where the boundary $\partial \omega= \{t \not\in \omega/ \exists s \in \omega, |t-s|<\ell \}$
and with 
$$\mathcal{N}=
\begin{cases}
N \mbox{ in the tetrode or non subcritical MEA case}\\
\kappa [\log(E/\alpha)]^2  \mbox{ in the  subcritical MEA case}
\end{cases}.$$

\item Moreover, if \begin{align}
\label{maxastarw}
\inf_{(n,t) \in \Sstar\cap \omega} |\vec{a}^*_{n,t}| > \frac{z_{\alpha} + \lambda + 2\|\vec{a}^*\|_{\infty,\partial \omega} \,(\rho \kl{\overline{c}} + \epsilon \mathcal{N})}{\cd - 2\epsilon  \,\mathcal{N}},
\end{align}
then $$
Supp(\ahat_\omega)= \Sstar\cap \omega.
$$
\end{enumerate}

\end{theorem}

See Section \ref{proof:lasso_results} for the proof.

\medskip

This theorem is stronger than the usual retrieval of support for Lasso estimator. Indeed it first applies to all possible subwindows at the same time, including the whole Lasso estimator itself. Next it does not calibrate $\lambda$ by the level of sparsity, that is $|\Sstar|$. Indeed in our problem even if $|\Sstar|$ is small compared to $T$, this grows linearly with $T$ since in expectation, under the assumptions of Proposition \ref{tempo}, it is roughly $pT$.

\medskip
Let us now discuss a bit more the choice of $\lambda$ and the calibration conditions. The first stringent condition is 
$$\cd>2\rho \pat{\overline{c}}+ 4 \epsilon \mathcal{N}.$$
Note that by Cauchy Schwarz , we already have that $\cd>\rho$, so what we ask here is a little bit stronger. The shape of the action potentials need to be picky enough to have $\rho$ small. In the same way, we need action potential shapes that are sufficiently different to have $\epsilon$ small.
The multiplication by $\mathcal{N}$ is in fact very large and a conservative upper-bound to the phenomenon taking place here. By looking at the proof, we can see that this is in fact the number of neurons (in a spatial overlap) that synchronizes with a lag less than $\ell$, that is  a few milliseconds. In practice, if this phenomenon is important for the neural coding \cite{tuleau2014multiple}, it usually involves a few neurons, except during epileptic crisis.

\medskip
Next, $\vec{a}^*$ is usually assumed to be a binary 0/1 vector in the classical problem. Therefore \eqref{minlambda} means  that we need
$$ \lambda > \frac{\cd+ 2\rho\pat{\overline{c}}}{\cd-2\rho\pat{\overline{c}} - 4\epsilon \mathcal{N}} \left( \sqrt{2\sigma^2 \co\log(2NT/\alpha)} +  2 (\rho\pat{\overline{c}} + \epsilon \mathcal{N}) \right),$$

On the other hand, Condition \eqref{maxastarw} becomes
$$\lambda<  \cd- 2\rho\pat{\overline{c}}-4\epsilon \mathcal{N} -  \sqrt{2\sigma^2 \co\log(2NT/\alpha)}.$$

So there is room to find such a $\lambda$ if typically 
$$\cd> \max\left(10\rho\pat{\overline{c}}+8\epsilon \mathcal{N}\quad , \quad  3\sqrt{2\sigma^2 \co\log(2NT/\alpha)}+4\rho\pat{\overline{c}}+6\epsilon \mathcal{N}\right).$$
This condition is reasonable in the context of spike sorting with high energy and peaky spike shapes, weak correlations between different neurons and normal neural activity.

\medskip
Moreover the windows generated by Algorithm \ref{algo:swas_detailed} are built to guarantee that $$\|\hat{\vec{a}}_\omega\|_{\infty,\partial \omega}=0,$$
which would imply that $$\|\vec{a}^*_\omega\|_{\infty,\partial \omega}=0$$
So in practice, $$ \lambda > \frac{\cd+ 2\rho\pat{\overline{c}}}{\cd-2\rho\pat{\overline{c}} - 4\epsilon \mathcal{N}}  \sqrt{2\sigma^2 \co\log(2NT/\alpha)}$$
should be sufficient.

\medskip
Also the windows which are generated by the algorithm can be controlled, as we can see in the following result.

\begin{cor}\label{control_omega}
Fix $\alpha\in (0,1/2)$.
Let Assumptions \ref{ref} and \ref{onG} be satisfied and let us assume that the noise variables are i.i.d. Gaussian. 
With the notations of Propositions \ref{spatial}, \ref{tempo} and Lemma \ref{conc},  on the same  event $\Omega_\alpha$ of probability larger than $1-\alpha$, if $\lambda$ is chosen so that \eqref{minlambda} and \eqref{maxastarw} are satisfied, then all the windows $\omega \in \Omega$ constructed by Algorithm  \ref{algo:swas_detailed} have a length controlled by 
$$\mathscr{W} =\begin{cases}
c' \log(T/\alpha) & \mbox{in the non-subcritical MEA or the tetrode case}\\ & \qquad \mbox{as soon as } p=Nm\Delta \leq 1/2,\\
c" \log(ET/\alpha) & \mbox{in the subcritical MEA  case}\\ & \qquad \mbox{as soon as the rate per cluster } p_c=N_c m\Delta \leq 1/2.
\end{cases}
$$
Also on the same event, steps 7,8 and 9 of the algorithm never occur.
\end{cor}

\subsection{Complexity of the sliding window working set algorithm}
\label{subsec:complexity}

We recall first an important existing result which gives the general (approximate) complexity of solving the Lasso with a working set algorithm. 
\begin{theorem}[\cite{loth2011}, Section 2.4]
\label{thm:general_complexity}
Consider the Lasso problem \eqref{eq:vec_model} with $n$ observations and $p$ features.
Then in order to compute a Lasso solution which selects $k$ features out of $p$, the working set algorithm has an approximate complexity of
    \begin{equation}
        \label{eq:general_complexity}
       C_{ws}(n,p,k) = O(n^2 p k + nk^3 + k^4).
    \end{equation}
\end{theorem}

Therefore applying Theorem~\ref{thm:general_complexity} to the resolution of the global Lasso problem with the working set strategy gives the following complexity.

\begin{proposition}
\label{prop:naive_activeset_complexity}
    Under the hypothesis of Theorem \ref{supp}, the algorithm \ref{algo:gas} solves a problem of $TE$ observations with $TN$ features and a number of true activations $O(TN)$, an application of Theorem \ref{thm:general_complexity} recovers an approximate
    complexity of
    \begin{equation}
        \label{eq:naive_complexity}
        C_{ws}(TE,TN,TN)= O(T^4(E^2N^2+EN^3+N^4)).
    \end{equation}
\end{proposition}

This result informs us that the naive global working set strategy cannot solve in
an efficient manner our problem. For a multi electrode array, the constants $E$
and $N$ are expected to be large. But even in the tetrode case, for which the
constants $E$ and $N$ remains small, the length of the signal $T$ might
arbitrarily increase depending on the duration of the experiment. The quartic complexity in $T$ and $N$ makes it impossible to apply this algorithm in practical situations.

\medskip
By contrast, we now state our result regarding the complexity of the sliding window working set.

\begin{theorem}
    \label{thm:complexity_swas}
    Under the hypothesis of Theorem \ref{supp} and Corollary \ref{control_omega}, there exists an event $\Omega_\alpha$ of probability larger than $1-\alpha$ such  that, on $\Omega_\alpha$,  the sliding window working set algorithm \ref{algo:swas_detailed} has the following approximate complexity 
    \begin{itemize}
        \item In the non-subcritical MEA case or in the tetrode case 
       \begin{equation}
            \label{eq:general_complexity_tetrode}
            O(T\log^4(T/\alpha) (E^2N^2+ E N^3 + N^4) ).
        \end{equation}
        \item In the subcritical MEA case  
    \end{itemize}
    \begin{equation}
        \label{eq:general_complexity_mea}
        O(ET\log^4(ET/\alpha) \log(E/\alpha)^8 ).
    \end{equation}

    \end{theorem}

    The above result reveals that our sliding window working set can avoid the high computational quartic costs from Proposition \ref{prop:naive_activeset_complexity} of the naive working set method thanks to the structure of the convolution. In addition our method achieves with high probability a very impressive quasi-linear time complexity $(T\log(T)^4)$ in $T$ for both tetrode and MEA. The complexity quadratic in $E$ and quartic in $N$ is still quartic in the non-subcritical MEA case or in the tetrode case but becomes quasi-linear with $O(E\log(E)^{12})$ in the  subcritical MEA case where the spatial overlaps limit the increase in size for the independent sub-problems.
    
    \medskip
    To the best of our knowledge, this is the first proof of complexity with high probability that recovers a quasi-linear complexity in the dimensionality of the data for solving the Lasso.

\section{Numerical experiments}
\label{sec:num_exp}

In this section we show with numerical simulations the performance of the
sliding window working set presented in section~\ref{sec:swas}. First we compare
its computation time with other approaches. As all these approaches
require to solve efficiently Lasso problems, we used the parallel implementation
of the accelerated proximal gradient FISTA \cite{beck2009fast} from
\cite{mairal2014spams}. Then we show the accuracy of the support of the Lasso
estimator
for several values of the regularization parameter and noise level, and also when the number of
synchronization between neurons grow. All experiments were performed of a simple notebook
having 8GB memory and a CPU
Intel(R) Core(TM) i7-4810MQ CPU @ 2.80GHz. The Python code from the experiments will be
made available on Github upon publication.

\subsection{Computational complexity}
\label{subsec:temporal_complexity}

In order to illustrate the performances of the sliding window working set, we
present here a comparison of the computation times of four different
approaches detailed below:

\begin{itemize}
	\item \textbf{Global solver}: This is a generic Lasso solver of
	\cite{mairal2014spams} using the accelerated proximal gradient FISTA to
	solve the global problem \eqref{eq:lasso_problem} with a
	pre-computed matrix $\tens{H}$. Since the size 
	of the design matrix $\tens{H}$ growths as $O(T^2)$, this approach rapidly
	suffers from the growth of $T$, both in terms of computation times and
	memory usage.
	\item \textbf{Working set (naive)}: This is a straightforward implementation of the working set algorithm~\ref{algo:gas} using FISTA as the inner Lasso solver. As previously stated, this approach does not scale properly in memory. Moreover, the computational complexity of the computation of the KKT is $O(ENT^2)$, which also does not scale well with the problem dimensionality.
	\item \textbf{Working set with convolution}: As discussed above, using
	standard solvers with a pre-computed matrix $\tens{H}$ is not scalable with
	the signal length $T$. In this method, we adapt the standard working set algorithm
	Algorithm~\ref{algo:gas} to take into account the
	structure of $\tens{H}$. 
	The computational bottleneck comes from
	the KKT violation (line 3 of Algorithm~{\ref{algo:gas}}). But in practice those
	conditions can be computed efficiently using a convolution with the shapes
	$\tens{W}$, instead of expensive matrix products. The residual
	$\vec{y}-\tens{H}\vec{\hat{a}}$ and the correlation $\tens{H}^T
	(\vec{y}-\tens{H}\vec{\hat{a}})$ (a convolution with reversed shapes) can be
	computed with complexity $O(ENT\ell)$ hence linear in $T$. The use of a working
	set also means that the storage of $\tens{H}$ is not necessary anymore since
	we solve the Lasso on the much smaller $\tens{H}_J$. We use the FISTA solver from
	\cite{mairal2014spams} to solve the sub-problems at each iteration. Finally
	note that $\tens{H}_J$ will be very sparse due to the convolutional model
	and the sub-problem can be solved on a matrix $\tens{\tilde H}_J$ of $O(E|J|\ell)$ lines instead of $O(ET)$.

										\item \textbf{Sliding window working set}: This is the method proposed in
	section~\ref{sec:swas} and described in Algorithm~\ref{algo:swas_detailed}.
	It focuses only on a small temporal window and slides the window when the
	problem is locally solved. Furthermore, since the research of the new activation is only carried out on the current window $\omega$ and not on the full signal, the computation cost of the KKT condition is greatly reduced from $O(ENT\ell)$ to $O(EN|\omega|\ell)$.
\end{itemize}

We have simulated our dataset realistically by using the classical model from
\cite{hodgkin1952} for the description of the shape of the action potentials,
and implemented by \cite{pouzathuxley}. In order to focus our study on the
influence of $T$, we limited ourselves to reasonable values for the number of
neurons ($N=5$) and the number of electrodes ($E=4$). Note that these small
sizes for the parameters would actually correspond to the resolution of the
problem on a single spatial group or to the tetrode case. 
\ld{We present in Table \ref{onG}, numerical values which correspond to the parameters introduced in Assumption \ref{onG}. These values depend a lot on the recordings in practice.}
\begin{table}[t]        
        \begin{center}        
            \begin{tabular}{|c|c|c|}            
                \hline
                Quantity & Mean & Stdev \\
                \hline
                $\epsilon$ & 173 & 175 \\
                \hline
                $\rho$ & 0.828 & 0.013 \\
                \hline
                $\underline{c}$ & 24.5 & 7.7 \\
                \hline
                $\overline{c}$ & 918 & 1052 \\
                \hline
            \end{tabular}
        \end{center}
        \caption{Orders of magnitude for the parameters introduced in Assumption \ref{onG}.}
\end{table}

\begin{figure}	\begin{center}
		\includegraphics[width=0.9\linewidth]{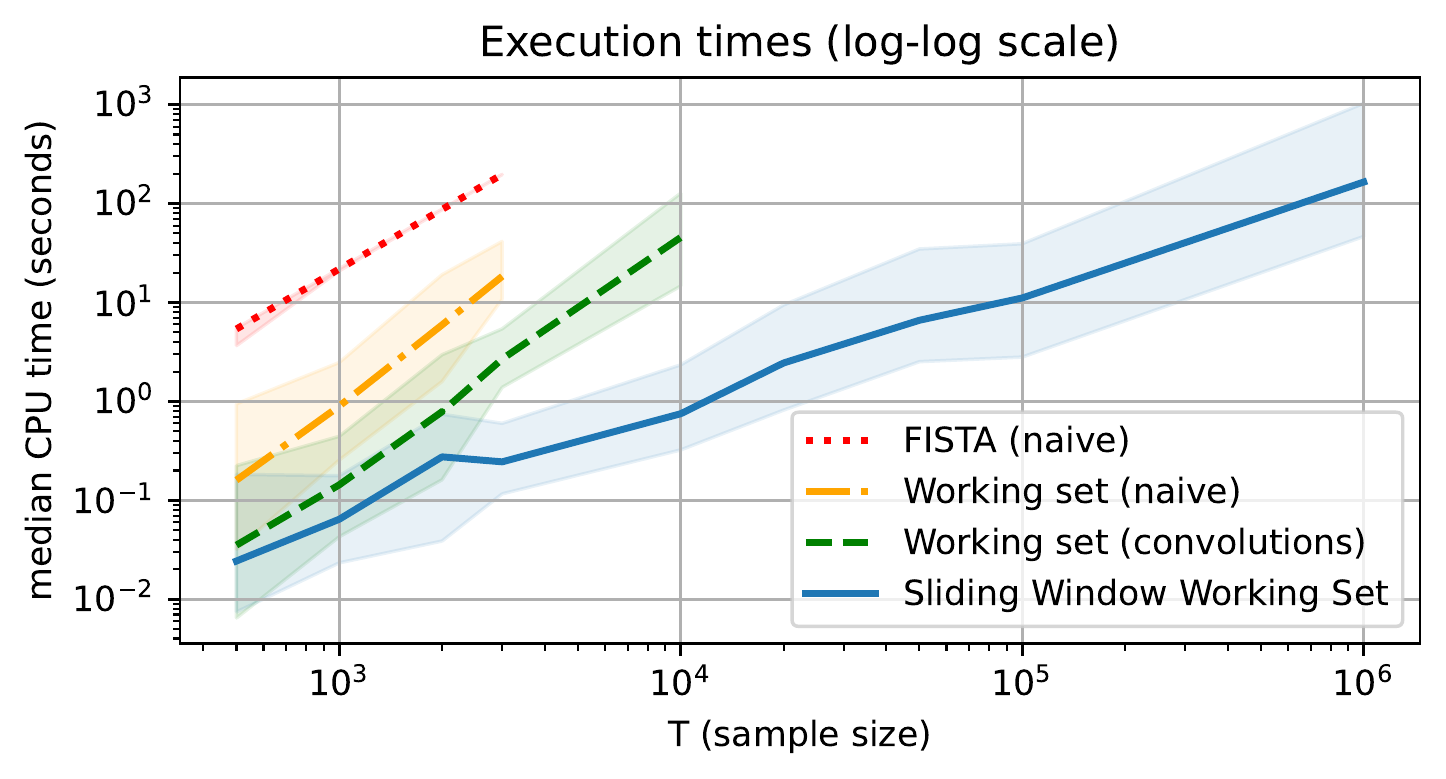}
		\caption{Comparison of the execution times for four different algorithms, when the size of the signals $n$ growths. We represented the median execution times over 40 simulations as the dotted lines. The bands represent the execution times between the 20 and 80 percentiles.}
		\label{fig:exec_times}
	\end{center}
\end{figure}

We present in Figure~\ref{fig:exec_times} the computation times of the different methods and their 20/80 percentiles  for different values of $T$ (each simulation is performed 40 times). It is clear from the figure that
the proposed algorithm is the most efficient and is actually the only one that
can solve problems with $T=10^6$ temporal samples. The slopes of the different
methods in the log-log space also show the differences in computational
complexities, with a slope near $1$ for the proposed algorithm. This corresponds to
the $O(T\log T)$ complexity obtained in the theoretical results.

Also note that we applied the algorithms to the same data, in order to verify that they all compute the same activation vectors, and therefore estimate the same supports.

\subsection{Influence of the noise and the regularization parameter}
\label{subsec:noise_lambda_influence}

Proper calibration of the regularization parameter $\lambda$ is crucial for the
success of the estimation. We want to visualize the influence of this choice,
especially for various noise levels. {Using real shapes of action potentials
recorded in  \cite{bethus2012} and that have been already spike sorted by
classical algorithms}, we simulate signals of size $T=500$ for different noise
levels, with $N=2$ neurons firing at $50Hz$ and recorded by $E=4$ electrodes. 

In order to measure the performance of the algorithm to recover the true support  we consider fist the classical F-measure used
for binary classification, which estimates a balance between false positive and
false negative rates. More precisely, we define the precision as $PRE =
\frac{TP}{TP+FP}$ and the recall as $REC = \frac{TP}{TP+FN}$, where $TP$, $FP$
and $FN$ are respectively the numbers of true positives, false positives and
false negatives. Then the F-measure is computed as $2 \frac{PRE.REC}{PRE+REC}$.
This measure tends to be pessimistic as it penalizes equally short and long
temporal deviations in the recovery. We introduce a softer measure of performance: $ CP(\vec{x},\vec{y})=
1 - \norm{K*(\vec{x}-\vec{y})}{1}/(\norm{\vec{x}}{1} + \norm{\vec{y}}{1})$,
where $K$ is a normalized rectangular function. Depending on the size of the
support of $K$, this measure allows us to penalize less small time deviations
than large time deviation. Here we took the size of its support equal to $10$, which corresponds to an accepted error of the order of 1ms.

\begin{figure}	\begin{center}
		\includegraphics[width=.8\linewidth]{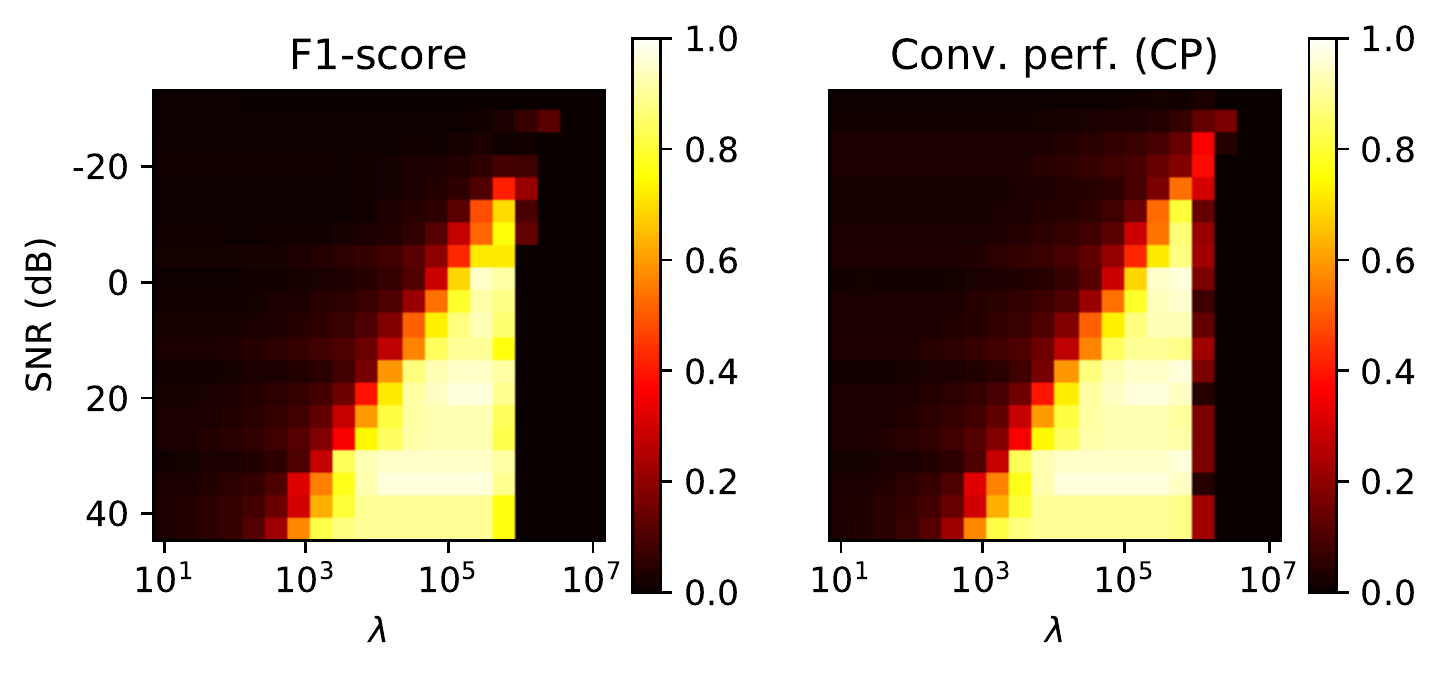}
		\caption{Influence of $\lambda$ and the signal-to-noise ratio on the performances of the Lasso estimator. Results are averaged over 5 draws.}
		\label{fig:noise_lambda_influence}
	\end{center}
\end{figure}

We provide on Figure~\ref{fig:noise_lambda_influence} the performance in F1
score (left) and the proposed convolutional performance (right). We can see that
the support recovery is very good in a large interval of values for the large
SNR but becomes narrow for low SNR where the support is harder to recover. High SNR recordings constitute an ideal setting for performing spike sorting, therefore the extracellular recording devices should be placed so that this SNR is high enough. Unfortunately this ideal environment may not always be guaranteed, especially in presence of bursting neurons, which action potential amplitudes may decrease down to the noise level \cite{lewicki1998review}. Therefore these experiments show that our method is robust enough to even treat low SNR recordings. \pat{Moreover even if the values in Table \ref{onG} are not very favorable to the method, we see that the method works well in practice.}

\subsection{Comparison with distance-based spike sorting methods}
\label{subsec:comparison_with_other_ss}

We now compare the performance of the Lasso estimator with distance based
methods that rely usually on K-means clustering for spike sorting. In a
clustering setting,the spike shapes $\tens{W}_n$ are the centroids of the method (\cite{ekanadham2014unified}). After the activation times have been estimated classical
approach select the neuron corresponding to the activation as the one closest to
the centroid. We now compare the Lasso and distance-based spike sorting in
the presence of synchronization between neurons (simultaneous spike). 

\begin{figure}
	\begin{center}
		\includegraphics[width=240pt]{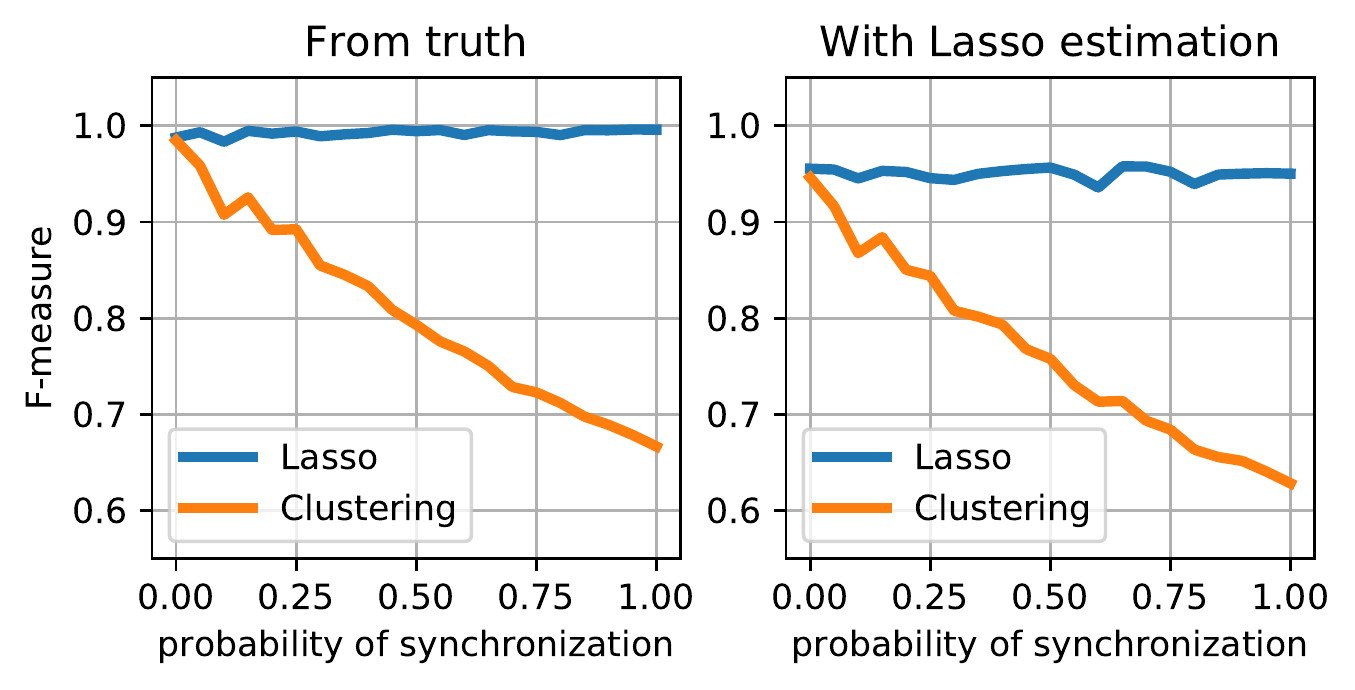}
		\caption{Comparison of Lasso and clustering performances (F-measure). Results are averaged over 50 draws. \label{fig:fmeasure_comparison}}
	\end{center}
\end{figure}

Using
similar data as in subsection~\ref{subsec:noise_lambda_influence}, we compare in
Figure \ref{fig:fmeasure_comparison} the performance of the  methods when the
number of synchronizations increases. To this end we use the true time
activation (left) where only which neurons are active is unknown and the support
recovered with our Lasso estimator The synchronizations have a minor impact on
the performances of the lasso estimator, illustrating the robustness of the
method due to the fact that the Lasso estimator is additive which mens that it
can handle simultaneous activation.

\section{Proofs}
\label{sec:proofs}

\subsection{Proof of Proposition \ref{spatial} \label{proof_spatial}}

In the sequel we use the term percolation term "cluster" to refer to spatial overlap.
  It has been shown in \cite{duminilcopin2018subcritical} that there exists a critical value $\gamma_c>0$ which only depends on $\delta$ and $r_0$, such that if $\gamma<\gamma_c$, the probability for a typical cluster to reach a radius $r$ (or a diameter $2r$) is less than $\exp(-c(\gamma,r_0)r)$, with $c(\gamma,r_0)>0$, depending on $\gamma$ and $r_0$.
  
  But the number of neurons, $N$, that can be sensed by the MEA is the number of neurons which are in a square of area $(\sqrt{E}+2)^2 \delta^2$. This is therefore a Poisson variable with mean $(\sqrt{E}+2)^2 \delta^2 \gamma$.
  
  So  by using basic concentration inequalities for Poisson variables (see for instance \cite{reynaud}), we obtain that, for all positive $x$
  $$\mathbb{P}(N> (\sqrt{E}+2)^2 \delta^2 \gamma+ (\sqrt{E}+2)\delta \sqrt{2 \gamma x}+ x/3)\leq e^{-x}.$$ 
  
  Let us take $e^{-x}=\alpha/3$ and let us enumerate the points (neurons) of the Poisson process in the whole plan with the first being the ones in the square of size $(\sqrt{E} +2) \delta$.
  We say that a cluster is attached to a neuron if the neuron belongs to this cluster.
  
  We use a union bound to control the size of each cluster attached to each neuron $n$ such that $n\leq Q$, with $Q$ the largest integer such that
  $$Q\leq  (\sqrt{E}+2)^2 \delta^2 \gamma+ (\sqrt{E}+2)\delta \sqrt{2 \gamma \log(3/\alpha)}+ \log(3/\alpha)/3. $$
  
  So we get that the probability to have one of these clusters of diameter larger than $r$ is  smaller than
  $$\left[(\sqrt{E}+2)^2 \delta^2 \gamma+ (\sqrt{E}+2)\delta \sqrt{2 \gamma \log(3/\alpha)}+ \log(3/\alpha)/3\right] e^{-c(\gamma, r_0) r}.$$
  
 We take  $r$ such that this bound is less than $\alpha/2$, that is 
 $$ r = \kappa' \log(E/\alpha),$$
 with
  $$\kappa'=\frac{1}{c(\gamma,r_0)} \left(1 + \frac{\log(2C)}{\log(\tau)}\right),$$
  and
  $$C = 4 \delta^2 \gamma+ 2 \delta \sqrt{\frac{2 \gamma}{\tau}}+ \frac{1}{3\tau}.$$
  which depend only on $\gamma,\delta,\tau$ and $r_0$.
  
  Finally, we can also control the number of neurons that belongs to each of the $Q$ balls that are used to encompass the $Q$ clusters of the first $Q$ points.
  
  With similar arguments as before on the control of Poisson variables and union bound, we can upper bound by $\alpha/3$, the probability that there is one of the $Q$ balls with more than $\kappa'' \log(E/\alpha)^2$ neurons in it.
  
  Therefore if we define $\Omega_{\alpha,s}$ as the event where (i) the total number of neurons is controlled (ii) the range of the $Q$ first clusters is controlled (iii) the number of neurons per ball for the first $Q$ balls is controlled, we obtain the desired result. 
 
\subsection{Proof of Proposition \ref{tempo} \label{proof_tempo}}

If $T_i$ is the $i$th index where $A_{T_i}=1$, then for all $i$, $\tau_i=T_i-T_{i-1}$ are independent Geometric variable on $\{1,2,...\}$ with parameter $p$, with the convention $T_0=0$.

We define $X_0=1$ and 
$$X_1=\min\{j>1,\tau_j>\eta\} \quad\mbox{and}\quad X_i=\min\{j>X_{i-1},\tau_j>\eta\}.$$
 Similarly, for $i\geq 1$, $\delta_i=X_i-X_{i-1}$ are independent geometric variables of parameter $(1-p)^\eta$.
 
 Therefore the $i$th overlap, which happens between $T_{X_{i-1}}$ and $T_{X_i}$ has a length $D_i= T_{X_i}-T_{X_{i-1}}+1$.
 
 So for all integer $k$,
 $$\mathbb{P}(D_i>k\eta+1)\leq \mathbb{P}(\delta_i>k) \leq (1-(1-p)^\eta)^k \leq (1-0.5^\eta)^k.$$
 
 We have at most $R$ overlaps with $R$ the largest integer such that $R\leq T/\eta$.
 
 By a union bound we can control all the $R$ first overlaps and the probability to have at least one overlap larger than $k\eta+1$ is controlled by
 $$T/\eta (1-0.5^\eta)^k.$$
 Forcing this last term to be $\alpha$ gives the value of $k$ and concludes the proof in the non subcritical case. The complementary event is $\Omega_{\alpha,A}$.
 
 In the subcritical case, using the notation of the proof of Proposition \ref{spatial}, we need to control it for all the first $Q$ clusters, which lead us to
 $$Q T_{max}/\eta (1-(1-p)^\eta)^k,$$
 hence the other choice of $k$. The complementary event is $\Omega_{\alpha,A}$.
 We then use here  $\Omega_{\alpha,A,s}=\Omega_{\alpha/2,s}\cap\Omega_{\alpha/2,A}$

\subsection{Proof of Lemma \ref{conc}}
\label{proof:lemmaconc}

For a fix $(n,t)$, since the noise is Gaussian, the random variable
$$
\vec{h}_{n,t}^\top \noise  =  \sum_{e=1}^E (\noise^e_t,...,\noise^e_{t+\l}) \pusht{\vec{w}_{n,e}}{t} 
$$
is also a Gaussian variable with variance bounded by $\overline{c} \sigma^2$. Thus the event 
$$\left|\vec{h}_{n,t}^\top \noise \right|\geq z$$
is of probability less than $2 e^{-z^2/(2\sigma^2 \co)}$.

\medskip

This argument is valid for any $(n,t)$. Therefore, by an union bound argument, we get, for a fixed $\alpha \in (0,1)$, with the choice  $z=z_\alpha$, for all $n$ in $1,..., N$ and all $t$ in $1,...,\T$,
 $$\left|\vec{h}_{n,t}^\top \noise \right|\leq z_\alpha,$$
with probability larger than $1-\alpha$. \qed

\subsection{Proof of Theorem \ref{supp}}
\label{proof:lasso_results} 

Note that if Algorithm \ref{algo:swas_detailed} needs to work with sparse matrices for computational reasons, mathematically speaking, we can as well work with the corresponding inflated matrices and this will not change the value of the solution, but just the space in which it is represented. Therefore, for the sake of convenience when we investigate the statistical properties of the method,
we define $\vec{a}_J$
as the vector obtained by setting to $0$ all the coordinates from $\vec{a}$ with their index not
in $J$. We define similarly $\tens{H}_J$ as the matrix obtained by replacing all the columns
from $\tens{H}$ with their index not in $J$ by the zero vector. Therefore in the sequel vectors $\vec{a}$ and matrix $\vec{H}$ have always the same dimensions.
\medskip

We work on the event $\Omega_\alpha$, which as stated in the remarks below Lemma \ref{conc} is of probability greater than $1-\alpha$. We fix a spatial overlap if we are in the subcritical MEA case or we work with the whole set of sensors in the other cases. In every cases, the number of neurons in the restricted problem is bounded by $\mathscr{N}$ thanks to Proposition \ref{spatial}. We also fix a given window $\omega$.

\medskip

We now solve the Lasso problem on the temporal window $\omega$ on the possibly restricted set of neurons:
$$
\hat{\vec{a}}_{\omega}=\arg\min_{\vec{a}/Supp(\vec{a})\subset \omega} \|\vec{y}-\tens{H}_{\omega}\, \vec{a}\|^2_2+ 2 \lambda \|\vec{a}\|_1,
$$
where we have used a slight abuse of language:  "$Supp(\vec{a})\subset \omega$" means that the temporal indices $t$ of $\vec{a}=(\vec{a}_{n,t})_{n,t}$ have to be in the temporal window $\omega$.
Let us define the solution of the Lasso optimization problem on $\Sstar \cap \omega$ where we recall that $\Sstar$ is the true support of $\vec{a}^*$:
$$
\hat{\vec{a}}_{\Sstar\cap \omega}
=\arg\min_{\vec{a}/Supp(\vec{a})\subset \Sstar\cap \omega} \|\vec{y}-\tens{H}_{\Sstar\cap \omega}\, \vec{a}\|^2_2+ 2 \lambda \|\vec{a}\|_1,
$$
where we have also made a slight abuse of language: $\Sstar \cap \omega =\{(n,t)\in \Sstar / t \in \omega\}$.
Note that $\hat{\vec{a}}_{\Sstar\cap \omega}$ and $\hat{\vec{a}}_\omega$ are both of dimension $NT$, with (temporal) support inside $\omega$.

\medskip

Our goal is to prove that $\hat{\vec{a}}_\omega$ is null outside of $\Sstar$. To this end, we first prove that the vector $\hat{\vec{a}}_{\Sstar\cap \omega}$ satisfies the KKT conditions on the temporal window $\omega$.

\medskip
Noting that
$
\H_{\omega}\,\hat{\vec{a}}_{\Sstar\cap \omega} = \H_{\Sstar \cap \omega}\,\hat{\vec{a}}_{\Sstar\cap \omega},
$ we see that $\hat{\vec{a}}_{\Sstar\cap \omega}$ already satisfies the KKT condition for any $(n,t) \in \Sstar\cap \omega$. We only have to check the KKT conditions for $(n,t) \in \omega\setminus \Sstar$ (with the same kind of language abuse as before).

\medskip

To this end, we first need to prove a bound on $\|\vec{a}^*-\hat{\vec{a}}_{\Sstar\cap \omega}\|
_{\infty,\omega}$ where $\|\vec{a}\|_{\infty,J} := \max_{n,t\in J} |a_{n,t}|$.

By definition, the Lasso solution $\hat{\vec{a}}_{\Sstar\cap \omega}$ satisfies the following necessary condition for all $(n,t)\in \Sstar\cap \omega$:
$$
|\vec{h}_{n,t}^\top (\vec{y}-\tens{H}_{\Sstar\cap \omega}\,\hat{\vec{a}}_{\Sstar\cap \omega})|\leq \lambda.
$$
We deduce that, for all $(n,t)\in \Sstar\cap \omega$,
$$
|\vec{h}_{n,t}^\top \tens{H}_{\Sstar\cap \omega}(\vec{a}^*- \hat{\vec{a}}_{\Sstar\cap \omega})|\leq \lambda + |\vec{h}_{n,t}^\top \noise| + |\vec{h}_{n,t}^\top \tens{H}_{\Sstar\cap \omega^c}\,\vec{a}^*|,
$$
with $\omega^c$ being the complementary in $\llbracket 1,T \rrbracket$ of the temporal window $\omega$.
\medskip

\medskip

In view of Lemma \ref{lem:Gblockband}, we have for all $(n,t)\in \Sstar\cap \omega$ that 
$$
\vec{h}_{n,t}^\top \tens{H}_{\Sstar\cap \omega^c}\,\vec{a}^* = \vec{h}_{n,t}^\top \tens{H}_{\Sstar\cap \partial\omega}\,\vec{a}^*. 
$$

In addition, Assumption \ref{onG} guarantees that $|G_{(n,t),(n',t')}|\leq \epsilon$ for all $n\neq n'$ and \ld{$|G_{(n,t),(n,t')}| \leq \rho |G_{(n,t),(n,t)}| \leq \rho \overline{c}$} for any $t\neq t'$. Also, given the refractory period, there can be at most only 1 activation on any interval of length $l$. Thus we get, for all  $(n,t)\in \Sstar\cap \omega$,
\begin{align}
\label{eq:partialomega}
|\vec{h}_{n,t}^\top \tens{H}_{\Sstar\cap \omega^c}\,\vec{a}^*|\leq 2 (\ld{\rho \overline{c}} + \epsilon \mathcal{N})\|\vec{a}^*\|_{\infty,\partial \omega}.
\end{align}

\medskip
\noindent
Next, we have for all $(n,t) \in \Sstar \cap \omega$
\begin{align*}
\vec{h}_{n,t}^\top\H_{\Sstar\cap\omega} \,( \vec{a}^*-\hat{\vec{a}}_{\Sstar\cap \omega}) &= \sum_{t'\,:\, (n,t')\in \Sstar\cap \omega} \tens{G}_{(n,t),(n,t')}(\vec{a}^*_{n,t'} - (\hat{\vec{a}}_{\Sstar\cap \omega})_{n,t'})     \\
&\hspace{1cm}+ \sum_{n'\neq n} \sum_{(n',t') \in  \Sstar\cap \omega} \tens{G}_{(n,t),(n',t')}(\vec{a}^*_{n',t'} - (\hat{\vec{a}}_{\Sstar\cap \omega})_{n',t'}).
\end{align*}
Given the block-band structure of the Gram matrix $G$ (see Lemma \ref{lem:Gblockband}), the first sum in the above display contains exactly 1 nonzero term corresponding to $t'=t$:
$$
\sum_{t'\,:\, (n,t')\in \Sstar\cap \omega} \tens{G}_{(n,t),(n,t')}(\vec{a}^*_{n,t'} - (\hat{\vec{a}}_{\Sstar\cap \omega})_{n,t'})  := \tens{G}_{(n,t),(n,t)}(\vec{a}^*_{n,t} - (\hat{\vec{a}}_{\Sstar\cap \omega})_{n,t} ).
$$
Regarding the second sum, we also have in view of Lemma \ref{lem:Gblockband}:
\begin{align*}
\sum_{n'\neq n} \sum_{(n',t') \in  \Sstar\cap \omega} \tens{G}_{(n,t),(n',t')}(\vec{a}^*_{n',t'} - (\hat{\vec{a}}_{\Sstar\cap \omega})_{n',t'}) &= \sum_{n'\neq n} \sum_{(n',t') \in  \Sstar\cap \omega, \, |t'-t|\leq l} \tens{G}_{(n,t),(n',t')}(\vec{a}^*_{n',t'} - (\hat{\vec{a}}_{\Sstar\cap \omega})_{n',t'})
\end{align*}
Set $\Delta_{n,t} := \vec{a}^\star_{n,t} - (\hat{\vec{a}}_{\Sstar\cap \omega})_{n,t} $. Combining the last four displays, we get for all $(n,t)\in \Sstar\cap \omega$,
\begin{align*}
    \tens{G}_{(n,t),(n,t)} |\Delta_{n,t}| \leq   \sum_{n'\neq n} \sum_{(n',t') \in  \Sstar\cap \omega, \, |t'-t|\leq l} |\tens{G}_{(n,t),(n',t')}|\,|\Delta_{n',t'}| + |\vec{h}_{n,t}^\top \noise| + \lambda + 2 (\ld{\rho \overline{c}} + \epsilon \mathcal{N})\|\vec{a}^*\|_{\infty,\partial \omega}.
\end{align*}

Assumption \ref{onG} guarantees that $|\tens{G}_{(n,t),(n',t')}|\leq \epsilon$ for all $n\neq n'$ and $G_{(n,t),(n,t)}>\cd$. Also for a given $n'$, because of the refractory period,  there is at most 2 activations for this particular neuron at distance $\ell$ of $t$. Thus we get, for all  $(n,t)\in \Sstar\cap \omega$,
\begin{align*}
    \cd \,|\Delta_{n,t}| \leq 2\epsilon \,\mathcal{N} \,\|\Delta\|_{\infty,\omega} +  \max_{(n,t)
\in \Sstar\cap \omega} |\vec{h}_{n,t}^\top \noise| + \lambda+ 2 (\ld{\rho \overline{c}} + \epsilon \mathcal{N})\|\vec{a}^*\vec{a}^*\|_{\infty,\partial \omega},
\end{align*}
and consequently, since we are on $\Omega_\alpha$,
\begin{align*}
    (\kl{\cd} - 2\epsilon  \,\mathcal{N}) \|\Delta\|_{\infty,\omega} \leq z_\alpha+ \lambda + 2 (\ld{\rho \overline{c}} + \epsilon \mathcal{N})\|\vec{a}^*\|_{\infty,\partial \omega}.
\end{align*}
Combining this result with Lemma \ref{conc}, we get\begin{align}
\label{eq:supDelta}
      \|\vec{a}^*- \hat{\vec{a}}_{\Sstar\cap \omega}\|_{\infty,\omega} \leq \frac{z_{\alpha} + \lambda + 2 (\ld{\rho \overline{c}} + \epsilon \mathcal{N})\|\astar\|_{\infty,\partial \omega}}{\cd - 2\epsilon  \,\mathcal{N}}.
\end{align}

\medskip
We now check the KKT conditions for $\hat{\vec{a}}_{\Sstar\cap \omega}$ on $\omega\setminus \Sstar$. For any $(n,t)\in \omega\setminus \Sstar$, we have
\begin{align*}
\vec{h}_{n,t}^\top (\vec{y} - \tens{H}_{\omega}\hat{\vec{a}}_{\Sstar\cap \omega}) &=  \vec{h}_{n,t}^\top (\tens{H} \vec{a}^* + \xi - \H_{\omega}\hat{\vec{a}}_{\Sstar\cap \omega})\\
 &=  \vec{h}_{n,t}^\top \, \H_{\omega}\,  (\vec{a}^* -\hat{\vec{a}}_{\Sstar\cap \omega}) + \vec{h}_{n,t}^\top \, \H_{\Sstar\cap \partial\omega}\, \vec{a}^* +\vec{h}_{n,t}^\top\, \noise.
\end{align*}
In view of (\ref{eq:partialomega}) and (\ref{eq:supDelta}), we have on the event $\Omega_{\alpha,A,s}\cap \Omega_{\alpha,\xi}$, for all $(n,t)\in \omega\setminus \Sstar$,
\begin{align}
\label{eq:KKTstrict1}
|\vec{h}_{n,t}^\top (\vec{y} - \H_{\omega}\hat{\vec{a}}_{\Sstar\cap \omega})| &\leq z_{\alpha} +  2 (\ld{\rho \overline{c}} + \epsilon \mathcal{N}) \left( \|\vec{a}^*\|_{\infty,\partial\omega}+   \frac{z_{\alpha} + \lambda + 2 (\ld{\rho \overline{c}} + \epsilon \mathcal{N})\|\vec{a}^*\|_{\infty,\partial \omega}}{\cd - 2\epsilon  \,\mathcal{N}}\right).
\end{align}
We need the following condition to satisfy the strict KKT conditions:
\begin{align}
\label{eq:KKTstrict2}
\lambda > z_{\alpha} +  2 (\ld{\rho \overline{c}} + \epsilon \mathcal{N})\left( \|\vec{a}^*\|_{\infty,\partial\omega}+   \frac{z_{\alpha} + \lambda + 2 (\ld{\rho \overline{c}} + \epsilon \mathcal{N})\|\vec{a}^*\|_{\infty,\partial \omega}}{\cd - 2\epsilon  \,\mathcal{N}}\right),
\end{align}
or equivalently
\begin{align*}
\lambda > \frac{1+\frac{2(\ld{\rho \overline{c}}+2\epsilon \mathcal{N})}{\cd-2\epsilon \mathcal{N}}}{1-\frac{2(\ld{\rho \overline{c}}+2\epsilon \mathcal{N})}{\cd-2\epsilon \mathcal{N}}} \left( z_{\alpha} +  2 (\ld{\rho \overline{c}} + \epsilon \mathcal{N}) \|\vec{a}^*\|_{\infty,\partial\omega} \right).
\end{align*}

This means that $\hat{\vec{a}}_{\Sstar\cap \omega}$ satisfies the strict KKT conditions in (\ref{eq:KKTstrict}) below on the temporal window $\omega$. Thus we proved, that $\hat{\vec{a}}_{\Sstar\cap \omega}$ is a solution of the Lasso minimization problem on the temporal window $\omega$, on the event $\Omega_\alpha$.

\medskip

The following property is an immediate consequence of the convexity of the Lasso objective function.

\begin{lemma}\label{classic}
Consider $Crit(\vec{a})= \| \vec{y}- \H_{\omega} \vec{a}\|^2_2+ 2 \lambda \|\vec{a}\|_1$. Let $\tilde{\vec{a}}_{\omega}$ be a minimizer of $Crit(\vec{a})$, hence satisfying the KKT conditions:
\begin{equation}
    \label{eq:KKTstrict}
    \begin{cases}
        \vec{h}_{n,t}^\top (\vec{y} - \H_{\omega} \tilde{\vec{a}}_{\omega}) = \lambda \ \text{sign}((\tilde{a}_{\omega})_{n,t}) & \text{, if } (\tilde{\vec{a}}_{\omega})_{n,t} \neq 0,\\
        | \vec{h}_{n,t}^\top (\vec{y} - \H_{\omega} \tilde{\vec{a}}_{\omega})| \leq  \lambda & \text{,if} \;(\tilde{\vec{a}}_{\omega})_{n,t} = 0.
    \end{cases}
\end{equation}
Let $$\tilde{S}=\{(n,t) / | \vec{h}_{n,t}^\top (\vec{y} - \H_{\omega}\tilde{\vec{a}}_{\omega})| = \lambda\}.$$

Then for any other minimizer $\hat{\vec{a}}_{\omega}$ of $Crit(\vec{a})$, we have
$$
Supp(\hat{\vec{a}}_{\omega}) \subset \tilde{S}.$$
\end{lemma}

\medskip

Note that $\tilde{S}$ might be larger than the true support of $\tilde{\vec{a}}_{\omega}$ because there might be coordinates $(n,t)$ such that $(\tilde{\vec{a}}_{\omega})_{n,t} = 0$  and for which  $| \vec{h}_{n,t}^\top (\vec{y} - \H_{\omega} \tilde{\vec{a}}_{\omega})| =  \lambda$.

\begin{proof}[Lemma \ref{classic}]

In view of (\ref{eq:KKTstrict}), we have for any $(n,t)\in \omega$,
$$
 \vec{h}_{n,t}^\top (\vec{y} - \H_{\omega} \tilde{\vec{a}}_{\omega}) = \lambda\, s_{n,t},
$$
where $(s_{n,t})_{(n,t)\in\omega}$ is such that,
$$
\begin{cases}
    |s_{n,t}|\leq 1  & \text{, in all cases,}\\
        s_{n,t}= \text{sign}((\tilde{a}_{\omega})_{n,t}) & \text{, if } (\tilde{\vec{a}}_{\omega})_{n,t} \neq 0,\\
        | s_{n,t}| <  1 & \text{, if } (n,t)\in \tilde{S}\kl{^c}.    \end{cases}
$$

Therefore, we have
\begin{eqnarray*}
Crit(\tilde{\vec{a}}_{\omega}+\vec{a})-Crit(\tilde{\vec{a}}_{\omega}) &=& \| \vec{y}- \H_{\omega} \,(\tilde{\vec{a}}_{\omega} +\vec{a}) \|^2_2 - \| \vec{y}- \H_{\omega} \,\tilde{\vec{a}}_{\omega} \|^2_2+ 2 \lambda (\|\tilde{\vec{a}}_{\omega} +\vec{a}\|_1-\|\tilde{\vec{a}}_{\omega}\|_1)\\
&=& \|\H_{\omega}\, \vec{a} \|^2_2 - 2 <\vec{y}-\H_{\omega} \tilde{\vec{a}}_{\omega}, \H_{\omega} \,\vec{a}>+ 2 \lambda (\|\tilde{\vec{a}}_{\omega} +\vec{a} \|_1-\|\tilde{\vec{a}}_{\omega}\|_1)\\
&=& \|\H_{\omega}\, \vec{a} \|^2_2 + 2 \lambda  ( \sum_{(n,t)\in \omega} |(\tilde{\vec{a}}_{\omega})_{n,t}+\vec{a}_{n,t}|-|(\tilde{\vec{a}}_{\omega})_{n,t}|-\vec{a}_{n,t} s_{n,t} ).
\end{eqnarray*}

Set $\vec{a}=\hat{\vec{a}}_{\omega} - \tilde{\vec{a}}_{\omega}$. By convexity of the $l_{1}$-norm, we have for all $(n,t)\in\omega$, 
$$
|(\tilde{\vec{a}}_{\omega})_{n,t}+\vec{a}_{n,t}|-|(\tilde{\vec{a}}_{\omega})_{n,t}| - s_{n,t} \vec{a}_{n,t} \geq 0.
$$
Thus,
\begin{align*}
\sum_{(n,t)} (|(\tilde{\vec{a}}_{\omega})_{n,t}+\vec{a}_{n,t}|-|(\tilde{\vec{a}}_{\omega})_{n,t}|- s_{n,t}  \vec{a}_{n,t} ) &\geq 0  \end{align*}

\medskip
Assume that $Supp(\hat{\vec{a}}_{\omega} )$ is not included in $\tilde{S}.$ Then there exists $(n_0,t_0)\in Supp(\hat{\vec{a}}_{\omega} ) \cap {\tilde{S}}^c$ such that $\vec{a}_{n_0,t_0} =(\hat{\vec{a}}_{\omega})_{n_0,t_0}\neq 0$ and $|s_{n_0,t_0}|<1$. Consequently, we get
$$
|\vec{a}_{n_0,t_0} | - s_{n_0,t_0} \vec{a}_{n_0,t_0}>0.
$$

\medskip
Since both $\hat{\vec{a}}_{\omega}$ and $\tilde{\vec{a}}_{\omega}$ are Lasso solution, we have
$$
0 = Crit(\hat{\vec{a}}_{\omega})-Crit(\tilde{\vec{a}}_{\omega})) \geq 2\lambda \, \sum_{(n,t)} |(\tilde{\vec{a}}_{\omega})_{n,t}+\vec{a}_{n,t}|-|(\tilde{\vec{a}}_{\omega})_{n,t}|-\vec{a}_{n,t} s_{n,t} ) \geq 2\lambda   (|\vec{a}_{n_0,t_0} | - s_{n_0,t_0} \vec{a}_{n_0,t_0})> 0.
$$
We obtain a contradiction. This means that 
$$
Supp(\hat{\vec{a}}_{\omega})\;\subset\;\tilde{S}. $$

\end{proof}

\medskip

\paragraph{Proof of Theorem \ref{supp} (continued)}

By Lemma \ref{classic} applied to $\tilde{\vec{a}}_{\omega}=\hat{\vec{a}}_{\Sstar\cap \omega}$, we get the first inclusion.

We assume in addition (\ref{maxastarw}). Then,
in view of (\ref{eq:supDelta}), we have on the event on $\Omega_\alpha$ $$
Supp(\hat{\vec{a}}_{\Sstar\cap \omega}) = \Sstar.
$$

\subsection{Proof of Corollary \ref{control_omega}}
With the choices provided in Algorithm \ref{algo:swas_detailed}, we start with the window $\omega=\llbracket 1,4\ell\rrbracket$. So $\omega$ is included in the first temporal overlap of size $\eta=4\ell$. Next the algorithm will compute and expand this window $\omega$ to the first time that no activation of $\hat{\vec{a}}_\omega$ is found in the last $2\ell$ coordinates. Thanks to Theorem \ref{supp}, this means that this is also the first time that $\vec{a}^*$ has no activation in a segment of length $2\ell$. This is not necessarily the first hole of size $2\ell$, because the algorithm only looks at $k\ell$ for some integer $k$, but definitely, the algorithm will stop and start a new window at the first "hole" of size $4\ell$.

In this sense, the first window will not be expanded not after the first hole of size $\eta$. Therefore its length is controlled by the control of the temporal overlap (see Proposition \ref{tempo}). The next step of the algorithm (Steps 7,8,9) cannot happen on the same event, indeed it would mean that the lasso estimator finds something at the beginning of the new window, whereas the estimator on the previous window (and therefore the truth) have no activation there. This is not possible since on every window, the Lasso estimator has the same support as the truth. 

Therefore we start a new window without merging with the one before and expand it again. The same arguments as before will apply recursively to prove our statement.

This conclude the proof of the theorem.

\subsection{Proof of complexity in Theorem \ref{thm:complexity_swas}}

In this proof we assume that the problem respects the  hypothesis of Theorem
\ref{supp} and Corollary \ref{control_omega}, there exists an event
$\Omega_\alpha$ of probability larger than $1-\alpha$. The following of the
proof suppose that we are on $\Omega_\alpha$.

\rf{We first investigate the non-subcritical MEA case or in the tetrode case.
Let us define $\tilde{T}=\log(T/\alpha)$ for the bounds on the window sizes.
In
this case we need to solve $O(T)$ temporal independent problems whose window size is
bounded by $O(\tilde{T})$ (as proven in Corollary \ref{control_omega}). Those independent problems, using notations from
Theorem \ref{thm:general_complexity}, have dimensionalities
of $O(E \tilde{T})$ observations, $O(N \tilde{T})$ features and
$O(N \tilde{T})$ selected features.
This means that the complexity is 
\begin{align*}
 \mathcal{C}(Alg. \ref{algo:swas_detailed}, Tetrode) & =  O(T C_{ws}(E \tilde{T},N \tilde{T},N \tilde{T})))\\
 &= O( T ((E \tilde{T})^2  (N \tilde{T})^2+ (E \tilde{T}) (N \tilde{T})^3+ (N \tilde{T}^4) ) )\\
 &= O(T\log(T/\alpha)^4 (E^2N^2+ EN^3 + N^4) )
  \end{align*}
which proves the result in equation \eqref{eq:general_complexity_tetrode}.

\medskip
In the subcritical MEA case, the problem can be solved with $O(ET)$ independent
problems (using both spatial and temporal overlaps). But those problems are of
much smaller size. In order to simplify the notations, let us now define $\tilde{T}=\log(ET/\alpha)$ and $\tilde{N}=\tilde{E}=\log(T/\alpha)^2$ for the bounds on the window sizes, $N_c$ and $E_c$.
Indeed Corollary \ref{control_omega} tells us that the size
of the temporal window is bounded in this case by $O(\tilde{T})$ and
Proposition \ref{spatial} tells us that the size of the spatial overlaps $N_c$
and $E_c$ are respectively bounded in this case by $O(\tilde{N})$ and $O(\tilde{E})$. This means that the
problems we need to solve have maximum dimensionality of
$O(\tilde{E} \tilde{T})$ observations,
$O(\tilde{N} \tilde{T})$ features and again
$O(\tilde{N} \tilde{T})$ selected features. This means that the
complexity of solving the whole problem is 
\begin{align*}
    \mathcal{C}(Alg. \ref{algo:swas_detailed}, MEA) & =  O(ETC_{ws}(\tilde{E} \tilde{T},\tilde{N} \tilde{T},\tilde{N} \tilde{T})) \\
     &= O(ET \tilde{T}^4 (\tilde{E}^2N^2+ \tilde{E} \tilde{N}^3 + \tilde{N}^4) ) \\
    &= O(ET\log(ET/\alpha)^4 \log(E/\alpha)^8 ).
      \end{align*}
   This proves result in equation \eqref{eq:general_complexity_mea}
   concludes the proof of the theorem.}

\section{Conclusion}
\label{sec:conclusion}

In this paper we propose a novel sliding window working set algorithm that can
solve exactly the large scale Lasso in spike sorting in an efficient way by
exploiting the convolutional structure of the problem. Under
some realistic assumptions, we prove that the complexity of the proposed algorithm
is quasi-linear with high probability with respect
to the temporal dimensionality of the signal.  Under some conditions on
the neurons firing rate, the complexity is also quasi-linear with the number of
electrodes which is a very important aspect for MEA that can have a large number
of electrodes.  We perform numerical experiments on a realistic signals and
recover the theoretical computational complexity.

We believe that this result opens the door for large scale spike sorting on the
recently available MEA sensor but also in other potential application which rely
on a sparse estimation with a convolutional model. Future work will investigate
the simultaneous estimation of the spike shape and activation and the online
update of those shapes along time.

\bibliographystyle{spmpsci}

\end{document}